\newif\ifdatestamp \newif\ifshowkeys 

\datestampfalse
\showkeysfalse

\documentclass[1p]{elsarticle}

\biboptions{sort} %This sorts a multiple reference to be [4,5,6] instead of [4,6,5].

\usepackage{hyperref}
\usepackage{thmref}
\usepackage{ahthesis}

\ifshowkeys\usepackage{showkeys}\fi

\usepackage{fancyhdr}
\newtheorem{theorem}{Theorem}[section]

\newtheorem{corollary}[theorem]{Corollary}
\newdefinition{definition}[theoremctr]{Definition}
\newtheorem{proposition}[theorem]{Proposition}

\newtheorem{lemma}[theorem]{Lemma}

\newproof{proof}{Proof}

\usepackage{amssymb}

\def\qed{\trueqed}

\begin{document}

\ifdatestamp
  \pagestyle{fancy}
\fi

%%%%%%%%%%%%%%%%%%%%%%%%%%%%
%% TITLE PAGE INFORMATION %%
%%%%%%%%%%%%%%%%%%%%%%%%%%%%

\thispagestyle{empty}

\title{Construction of supercharacter theories of finite groups}
\author{Anders O.F.~Hendrickson}
  \ead{ahendric@cord.edu}
  \address{Department of Mathematics and Computer Science,
           Concordia College, 
           %901 8th St S, 
           Moorhead, MN 56562, USA}

\begin{abstract}

  Much can be learned about a finite group from its character table,
  but sometimes that table can be difficult to compute.
  Supercharacter theories are generalizations of character theory,
  defined by P.~Diaconis and I.M.~Isaacs in \cite{diaconis_isaacs},
  in which certain (possibly reducible) characters called supercharacters
  take the place of the irreducible characters,
  and a certain coarser partition of the group
  takes the place of the conjugacy classes.
  
  In particular, if $\KK$ is a partition of a finite group $G$, 
  there may exist a compatible partition $\XX$ of the irreducible characters of $G$, 
  along with a character $\chi_X$ for every $X\in\XX$ with the elements of $X$ as its irreducible constituents, 
  so that each $\chi_X$
  is constant on each $K\in\KK$ and $|\XX| = |\KK|$. 
  If every irreducible character is a constituent of some $\chi_X$, 
  then the ordered pair
  $(\XX,\KK)$ is called a supercharacter theory.

  We present five new ways to construct new supercharacter theories
  out of supercharacter theories already known to exist,
  including a direct product, 
  a lattice-theoretic join,
  two products over normal subgroups, 
  and a duality for supercharacter theories of abelian groups.
\end{abstract}


\maketitle

\ifdatestamp
  \pagestyle{fancy}
\fi

%\chapter{1}{Introduction}\label{ch_introduction}

\section{Introduction}

Let $G$ be a finite group with $n$ conjugacy classes,
and let $\Irr(G)$ be the set of complex irreducible characters of $G$.
Now the $n$ conjugacy classes partition the group $G$,
and each of the $n$ irreducible characters
is constant on every conjugacy class of $G$.
The relationship between the irreducible characters
and the conjugacy classes is of great use in studying a finite group,
but for some groups $\Irr(G)$ is quite difficult to compute.

For example, let $U_n(\F_q)$ be the group
of upper triangular matrices over the field of size $q$
with all diagonal entries one.
The irreducible characters of $U_n(\F_q)$ are quite difficult to compute,
but
Carlos Andr\'e \cite{andre1995, andre1999, andre2002}
and Ning Yan \cite{yan}
developed 
a theory of ``basic characters'' or ``transition characters'' 
as an approximation to the full character table of $U_n(\F_q)$.
In this theory certain reducible characters take the place of irreducible characters
and the role of conjugacy classes is played by certain unions of conjugacy classes.
Although these theories are simple enough to be computed explicitly,
%These theories resemble the ordinary theory of irreducible characters
%in some ways 
they are rich enough to handle some problems
that traditionally required knowing the full character theory
\cite{arias2004}.
Persi Diaconis and I.~Martin Isaacs \cite{diaconis_isaacs} 
have generalized the work of Andr\'e and Yan
to the concept of a supercharacter theory of a finite group,
which is defined as follows.
All groups mentioned will be finite.

%\pagebreak
\begin{definition}\label{defn_supercharthy}
  Let $G$ be a finite group,
  let $\KK$ be a partition of $G$,
  and let $\XX$ be a partition of the set $\Irr(G)$.
  Suppose that for every part $X\in\XX$ there exists a character $\chi_{_X}$
  whose irreducible constituents lie in $X$,
  and suppose the following three conditions hold.
  \begin{enumerate}
    \item Each of the characters $\chi_{_X}$
          is constant on every part $K\in\KK$.
    \item $|\XX|=|\KK|$.
    \item Every irreducible character is a constituent of some $\chi_X$.
%    \item The set $\{1\}$ is one of the parts of $\KK$.
  \end{enumerate}
  Then 
  we call the characters $\chi_{_X}$ \defnstyle{supercharacters},
  we call the members of $\KK$ \defnstyle{superclasses},
  and we say that the ordered pair $(\XX,\KK)$ 
  is a \defnstyle{supercharacter theory}.
  If $\CCC=(\XX,\KK)$ is a supercharacter theory,
  then we define $|\CCC|$ to be the integer equal to
  both $|\XX|$ and $|\KK|$.
  We write $\Sup(G)$ for the set of all supercharacter theories of $G$.
\end{definition}

Assume that $(\XX,\KK)$ is a supercharacter theory of a group $G$,
and for every subset $X$ of $\Irr(G)$ let $\s_X$ be the character $\sum_{\psi\in X}\psi(1)\psi$.
Diaconis and Isaacs
prove in \cite[Lemma 2.1]{diaconis_isaacs}
that $\{1\}\in\KK$, that $\{1_G\}\in\XX$,
and that for each $X\in\XX$, 
the supercharacter $\chi_{_X}$ must be a constant multiple of $\s_X$.
It is therefore no loss to assume that $\chi_{_X}=\s_X$,
and we shall make that assumption throughout this paper.
It is also shown in 
\cite[Theorem 2.2(c)]{diaconis_isaacs}
that if $\CCC=(\XX,\KK)$ is a supercharacter theory,
then each of the partitions $\XX$ and $\KK$
uniquely determines the other.

This paper investigates several ways to discover
new supercharacter theories from supercharacter theories
already known to exist,
as an aid to future work classifying $\Sup(G)$ for given groups $G$.
Let us begin by reviewing the five constructions of supercharacter theories
given in \cite{diaconis_isaacs}.
First, there are two trivial supercharacter theories.
\begin{definition}\label{defn_minmaxthy}
  Let $G$ be a group.  
  Then the \defnstyle{minimal supercharacter theory} $\m(G)$
  is given by the partitions
  of $\Irr(G)$ 
  into singleton sets
  and of $G$ into its conjugacy classes.
  The \defnstyle{maximal supercharacter theory} $\M(G)$,
  on the other hand,
  is given by the partition $\{\{1_G\},\, \Irr(G)-\{1_G\}\}$ of $\Irr(G)$
  and the partition $\{\{1\},\,G-\{1\}\}$ of $G$.
\end{definition}

Because the superclasses of the minimal theory 
are the conjugacy classes
and its supercharacters are scalar multiples of the irreducible characters of $G$,
the minimal theory is just the ordinary character theory of $G$.
The supercharacters of the maximal theory, by contrast,
are the principal character $1_G$ and $\rho_G-1_G$,
where $\rho_G$ is the regular character of $G$.
It is routine to verify that these two theories satisfy 
the conditions of 
Definition \ref{defn_supercharthy}.
%\thmref{defn_supercharthy}.

Next, let $A$ be a group acting on $G$ by automorphisms.
Then $A$ permutes both $\Irr(G)$ and the set of conjugacy classes of $G$.
The partition of the conjugacy classes into $A$-orbits
yields a partition $\KK$ of the group $G$,
and taking $\XX$ to be the partition of $\Irr(G)$ into $A$-orbits,
a straightforward calculation shows
that for each $A$-orbit $X\in\XX$,
the character $\s_X$ is constant on each member of $\KK$.
A lemma of Richard Brauer guarantees that $|\XX|=|\KK|$,
and %the set $\{1\}$ lies in $\KK$ because automorphisms fix the group identity.
therefore $(\XX,\KK)$ is a supercharacter theory of $G$.

For the fourth construction given in \cite{diaconis_isaacs},
let $G$ be a group and let $H$ be a group of field automorphisms of $\C$;
then a supercharacter theory can be obtained
by taking $\XX$ to be the orbit decomposition of $\Irr(G)$ 
under the action of $H$.
Finally, 
the bulk of \cite{diaconis_isaacs}, like the subsequent papers by 
Diaconis, Otto, Thiem, Venkateswaran, and Marberg
\cite{thiem_diaconis, ottothesis, thiem_marberg, thiem_venkateswaran},
generalizes the theory of Andr\'e and Yan
to a particularly well-behaved supercharacter theory
of a certain family of groups called algebra groups.
In other papers, Andr\'e and Neto have recently begun to describe
a supercharacter theory for Sylow subgroups of symplectic and orthogonal groups
\cite{andreneto2006, andreneto2008a, andreneto2008b}.

In this article we present several new ways 
to obtain supercharacter theories
of a finite group $G$.
Section \ref{sect_directproducts} defines direct products of supercharacter theories,
and Sections \ref{sect_subalg} and \ref{sect_joins} 
derive a lattice-theoretic join operation on supercharacter theories.
Sections \ref{sect_starproduct} through \ref{sect_wtp} introduce
new supercharacter theories of a group $G$ as products
of supercharacter theories of quotient groups and normal subgroups of $G$,
and investigate properties of these products.
Restricting our attention to abelian groups,
in the closing sections we describe a duality relation
between the supercharacter theories of an abelian group $G$
and those of the group $\Irr(G)$ of its irreducible characters.

\ifdatestamp
  \pagestyle{fancy}
\fi

%\chapter{2}{Basic results}\label{ch_basicresults}

\section{Direct products}\label{sect_directproducts}

New supercharacter theories can often be found
by combining supercharacter theories that are already known.
Our first approach will be to form a direct product of supercharacter theories.
Let $M$ and $N$ be groups, and let $G$ be the direct product $M\by N$;
then we know that $\Irr(G) = \Irr(M)\by\Irr(N)$.
Given two supercharacter theories $(\XX,\KK)\in\Sup(M)$ and $(\YY,\LL)\in\Sup(N)$,
we can form a natural ``product'' theory $(\XX,\KK)\by(\YY,\LL)$
by forming a new supercharacter for every element of 
the cartesian product $\XX\by\YY$
and a new superclass for every element of $\KK\by\LL$.  Let
%$$\ZZ %= \XX\by\YY 
%      = \{ \{ \phi\by\theta : \phi\in X, \theta\in Y\} : X\in\XX, Y\in\YY\}$$
$$\ZZ = \{ X\by Y : X\in\XX,~Y\in\YY\} \mbox{~where~}
        X\by Y = \{ \phi\by\theta : \phi\in X,~\theta\in Y\}\sseq\Irr(G) $$
and let
%$$\MM %= \KK\by\LL 
%      = \{ \{ (m,n) : m\in K, n\in L\} : K\in\KK, L\in\LL \}.$$
$$\MM = \{ K\by L : K\in\KK,~L\in\LL \} \mbox{~where~}
       K\by L = \{ (m,n) : m\in K,~n\in L\}\sseq G.$$

\begin{lemma}
Using the above notation, $(\ZZ,\MM)\in\Sup(G)$.
\end{lemma}
\begin{proof}
Certainly $|\ZZ| = |\XX||\YY| = |\KK| |\LL| = |\MM|$.
%and since $\{1\}\in\KK$ and $\{1\}\in\LL$,
%it follows that $\{1\}\in\MM$.
%
So it suffices to show for all sets $X\in\XX$, $Y\in\YY$, $K\in\KK$, and $L\in\LL$
that the character $\sigma_{X\by Y}$ 
is constant on the set $K\by L$.
For all $m\in M$ and $n\in N$, we have
\begin{eqnarray*}
%\sigma_{XY}(g) &=& 
\sigma_{X\by Y}((m,n))
    &=& \sum_{\phi\in X} \sum_{\theta\in Y} (\phi\by\theta)((1,1))\cdot (\phi\by\theta)((m,n)) \\
    &=& \sum_{\phi\in X} \sum_{\theta\in Y} \phi(1)\theta(1)\phi(m)\theta(n) \\
%    &=& \sum_{\phi\in X} \sum_{\theta\in Y} \phi(1)\phi(m)\theta(1)\theta(n) \\
    &=& \sum_{\phi\in X} \phi(1)\phi(m) \sum_{\theta\in Y} \theta(1)\theta(n) \\
    &=& \sigma_X(m) \sigma_Y(n).% \\
    %&=& (\sigma_X \by \sigma_Y) ((m,n)).
    %&=& (\sigma_X \by \sigma_Y) (g).
\end{eqnarray*}
Thus for all $m,m'\in K$ and all $n,n'\in L$,
$$\sigma_{X\by Y}((m,n))
   = \sigma_X(m)\sigma_Y(n)
   = \sigma_X(m')\sigma_Y(n')
   = \sigma_{X\by Y}((m',n')).$$
Thus $\sigma_{X\by Y}$ is in fact constant on $K\by L$.
We conclude that $(\ZZ,\KK)$ is indeed a supercharacter theory of $G$.\qed
\end{proof}

Of course, not every supercharacter theory of a direct product
%\marginpar{This \P\ optional.}
arises in this manner; in particular,
the maximal supercharacter theory $\M(M\by N)$ does not lie in $\Sup(M)\by\Sup(N)$
unless either $M$ or $N$ is trivial.

\def\BB{\mathcal{A}}
\def\Ab{A}
\def\AA{A}

\section{Supercharacter theories and subalgebras}\label{sect_subalg}

%Given a group $G$, we would like to understand $\Sup(G)$ fully.
%\marginpar{Define $\sim$, $[~]$, and $[~]_\CCC$ before this section.}

To prepare for our next construction,
we investigate the connection between supercharacter theories of a finite group $G$
and certain subalgebras of the center of the group algebra $\C[G]$.
For every subset $K\sseq G$,
let $\Khat=\sum_{x\in K} x$;\index{$\Khat$}
recall that for every subset $X\sseq\Irr(G)$,
$\sigma_X$ denotes $\sum_{\psi\in X}\psi(1)\psi$.\index{$\sigma_X$}
Also recall that every character $\chi\in\Irr(G)$
has a corresponding central idempotent
$e_\chi = \rec{|G|}\chi(1)\sum_{g\in G} \overline{\chi(g)}\,g$,\index{$e_\chi$}
and that the set $\{e_\chi: \chi\in\Irr(G)\}$ is a basis for $\Zb(\C[G])$.
For every subset $X\sseq\Irr(G)$,
let $f_X = \sum_{\psi\in X} e_\psi$.\index{$f_X$}
Then $f_X$ and $\s_X$ are closely related.
Because $e_\chi = \rec{|G|}\sum_{g\in G}\overline{\s_{\{\chi\}}(g)}\,g$,
by the linearity of the $\sigma$ operator we have
\begin{equation}\label{eqn_fXsXreln}
f_X   = \summ_{\chi\in X} e_\chi
      = \rec{|G|} \sumg \bar{\s_X(g)} \, g .
\end{equation}

Diaconis and Isaacs show in \cite[Theorem 2.2(b)]{diaconis_isaacs} that if
$G$ is a group and if $\CCC=(\XX,\KK)\in\Sup(G)$,
then the set of superclass sums $\{\Khat: K\in\KK\}$
and the set of sums of idempotents $\{f_X: X\in\XX\}$
are two bases for the same subalgebra of $\Zb(\C[G])$,
which we shall denote $\BB(\CCC)$.

\label{sect_centralalgsuffices}
For each partition $\XX$ of $\Irr(G)$, 
%we can 
let $\Ab_\XX$ denote the subspace $\spanof{f_X: X\in\XX}$.  
Because the $f_X$ are central idempotents and
$f_X f_Y = 0$ if $X\neq Y$,
the subspace $\Ab_\XX$ is actually
a subalgebra of $\Zb(\C[G])$.
For example, 
if $\CCC=(\XX,\KK)$ is a supercharacter theory of $G$,
then $\Ab_\XX=\BB(\CCC)$.
It turns out that every subalgebra of $\Zb(\C[G])$ arises 
from a partition of $\Irr(G)$ in this way.

\begin{lemma}\label{lem_subalgyieldsfZs}
  %Phrase as bijection?
  Let $G$ be a group and let $\AA$ be a subalgebra of $\Zb(\C[G])$.
  Then there exists a unique partition $\ZZ$ of $\Irr(G)$
  such that $\{f_Z: Z\in\ZZ\}$ is a basis for $\AA$.
\end{lemma}
\begin{proof}
Because $\Zb(\C[G])$ is isomorphic to a direct sum of copies of $\C$,
it contains no nilpotent elements,
so neither does its subalgebra $\AA$;
hence the Jacobson radical $\Jb(\AA)=0$.
Then by Wedderburn's theorem
$\AA$ is a direct sum of full matrix rings;
but since $\AA$ is commutative, 
those are rings of $1\by 1$ matrices,
so $\AA$ too is a direct sum of copies of $\C$.
Hence $\AA$ is the linear span of some idempotents $f_1,\ldots,f_r$
whose sum is 1 and whose pairwise products are 0.
But every idempotent in $\Zb(\C[G])$ is a sum of some distinct $e_\chi$,
and because 
%$\sum_{i=1}^r f_i = 1 = \sum_{\mathclap{\chi\in\Irr(G)}} e_\chi$ 
$\sum_{i=1}^r f_i = 1 = \sum_{{\chi\in\Irr(G)}} e_\chi$ 
but the product $f_i f_j=0$ for $i\neq j$, 
every $e_\chi$ must appear in exactly one $f_i$.
Then there exists a partition $\ZZ$ such that 
$\{f_Z: Z\in\ZZ\}=\{f_1,\ldots,f_r\}$,
and this is the desired basis for $\AA$.

To show uniqueness, suppose $\YY$ is also a partition of $\Irr(G)$
such that 
$$\spanof{f_Y: Y\in\YY}=\AA=\spanof{f_Z:Z\in\ZZ}.$$
Let $\chi\in\Irr(G)$, and let $Y_0\in\YY$ and $Z_0\in\ZZ$ 
be the parts containing $\chi$.
Then because $f_{Y_0}\in\AA=\spanof{f_Z:Z\in\ZZ}$, 
the set $Y_0$ must be a union of parts of $\ZZ$;
in particular $Z_0\sseq Y_0$.
But by symmetry $Y_0\sseq Z_0$, so $Y_0=Z_0$.
Since the parts of $\YY$ and $\ZZ$ containing $\chi$ are identical
for all $\chi\in\Irr(G)$, it follows that $\YY=\ZZ$.\qed
\end{proof}

This \thmtype{lem_subalgyieldsfZs} enables us to determine whether an arbitrary partition $\KK$ of $G$
corresponds to a supercharacter theory,
using only computations in $\C[G]$
and making no mention of characters.

\begin{proposition}\label{prop_centralalgsuffices}
  Let $\KK$ be a partition of $G$. % with $\{1\}\in\KK$.
  Then there exists a partition $\XX$ of $\Irr(G)$ such that $(\XX,\KK)\in\Sup(G)$
  if and only if $\textspanof{\Khat: K\in\KK}$ is a subalgebra of $\Zb(\C[G])$.
\end{proposition}
\begin{proof}
As noted above, if $(\XX,\KK)$ is a supercharacter theory, then Diaconis and Isaacs
have proved \cite[Theorem 2.2(b)]{diaconis_isaacs}
that $\textspanof{\Khat: K\in\KK}$ is a subalgebra of $\Zb(\C[G])$.
So now suppose $\textspanof{\Khat: K\in\KK}$ is a subalgebra $\AA$ of $\Zb(\C[G])$.
%Note that $\{1\}$ must be a part of $\KK$ since $1\in\AA$.
Because $\AA\sseq\Zb(\C[G])$, each part $K\in\KK$ is a union of conjugacy classes of $G$.
Because the parts of $\KK$ are disjoint, we know that 
the set $\{\Khat: K\in\KK\}$ is linearly independent, and hence a basis for $\AA$.
By \thmref{lem_subalgyieldsfZs}, there exists some
partition $\XX$ of $\Irr(G)$ so that $\{f_X: X\in\XX\}$ is also a basis for $\AA$.
So $|\XX|=\dim\AA = |\KK|$, and it only remains to show that $\s_X$ is constant on $K$
for all $X\in\XX$ and all $K\in\KK$.
Now an element $y\in\C[G]$ lies in $\textspanof{\Khat:K\in\KK}$ if and only if
the function from $G$ to $\C$ which maps each group element to its coefficient in $y$
is constant on all $K\in\KK$.
But 
$$f_X = \sumg \left(\rec{|G|} \bar{\s_X(g)}\right) g$$
as we saw in Eq.~(\ref{eqn_fXsXreln});
then because $f_X\in\AA$, 
the function $g\mapsto \rec{|G|}\bar{\s_X(g)}$ 
is constant on all $K\in\KK$,
so $\s_X$ must be constant on all $K\in\KK$ as well.
Therefore $(\XX,\KK)$ is a supercharacter theory of $G$.\qed
\end{proof}

\section{Joins of supercharacter theories}\label{sect_joins}

Our next method of obtaining new supercharacter theories
is to perform a lattice-theoretic ``join'' on
two supercharacter theories already known.
Let us
first recall some well-known facts about partitions of a set.

For each finite set $S$,
the set $\Part(S)$ of partitions of $S$ into disjoint subsets
forms a partially ordered set under the relation ``$\pleq$''
where $\XX\pleq\YY$ if and only if every part of $\XX$ 
is contained in some part of $\YY$.
This poset in fact forms a lattice, 
which is called the \defnstyle{partition lattice} of $S$.
(See \cite[pp. 192ff]{gratzer1978} for more details.)
Thus for any two partitions $\XX$ and $\YY$ of $S$,
their join $\XX\join\YY$ and meet $\XX\meet\YY$ are defined.
The statement $\XX\pleq\YY$ is equivalent 
          to the statement $\XX\join\YY=\YY$
          and to the statement $\XX\meet\YY=\XX$.
We also note without proof the following easy facts.

\begin{lemma}\label{sorites}
  Let $\XX,\YY\in\Part(S)$.  Then the following hold:
  \begin{enumerate}
%    \item $\XX\pleq\YY$ if and only if $\XX\join\YY=\YY$,
%          if and only if $\XX\meet\YY=\XX$.\label{sorites0}
%    \item We have $[a]_\XX \sseq [a]_{\XX\join\YY}$ for all $a\in S$.\label{sorites1}%
%                \marginpar{Do we use \#\ref{sorites1}?}
    \item Suppose $\XX\pleq\YY$.
          Then each part of $\YY$ is the union of some parts of $\XX$.\label{sorites2}
          %\marginpar{Do we use \#\ref{sorites2}? Yes}
%    \item $a\sim_{\XX\join\YY} b$ iff there exist $g_0,\ldots,g_r\in S$
%          such that 
%          $$a=g_0\simx g_1\simy g_2 \simx \cdots \simx g_{r-1}\simy g_r = b.$$
%          \label{sorites_chain}
  \item Let $T$ be a set and let $f: S \to T$ be a function.
        Suppose $f$ is constant on each member of $\XX$ and each member of $\YY$.
        Then $f$ is constant on each member of $\XX\join\YY$.\label{fconstonxcupy}
  \end{enumerate}
\end{lemma}
\gobble{%OMIT PROOF
                     \begin{proof}
                     %The proof of part (\ref{sorites0}) is routine.
                     %
                     Part (\ref{sorites2}) holds since 
                     $Y\sseq \bigcup_{a\in Y} [a]_\XX \sseq Y$ for all $Y\in\YY$.
                     For part (\ref{fconstonxcupy}),
                     define an equivalence relation by saying
                     $a\sim b$ if and only if $f(a)=f(b)$,
                     and let $\KK$ be the partition of $S$ into the equivalence classes.
                     Then $f$ is constant on each member of a partition $\LL$
                     if and only if $\LL\pleq \KK$.
                     Then our hypotheses tell us that $\XX,\YY\pleq \KK$,
                     so their join $\XX\join\YY$ is also $\pleq \KK$.
                     Therefore $f$ is constant on each member of $\XX\join\YY$.\qed
                     \end{proof}
}%END GOBBLE

We saw in Section \ref{sect_centralalgsuffices}
that to every partition $\XX$ of $\Irr(G)$
corresponds a central subalgebra $\Ab_\XX=\textspanof{f_X: X\in\XX}$,
and that every central subalgebra arises in this way.
This bijection 
between $\Part(\Irr(G))$ and the subalgebras of $\Zb(\C[G])$
is also order-reversing with respect to partial orders.

\begin{lemma}\label{lem_fGalois}
The map $\XX\mapsto \Ab_\XX$ is a bijection 
from $\Part(\Irr(G))$ to the set of subalgebras of $\Zb(\C[G])$.
Moreover,
if $\XX,\YY\in\Part(\Irr(G))$,
then
\begin{enumerate}
\item $\XX\pleq\YY$ if and only if $\Ab_\YY\sseq \Ab_\XX$. \label{part_fGalois_po}
\item $\Ab_{\XX\join\YY} = \Ab_\XX \cap \Ab_\YY$. \label{part_fGalois_XcupY}
\item $\Ab_{\XX\meet\YY} = \gen{\Ab_\XX, \Ab_\YY}$,
      the subalgebra generated by $\Ab_\XX$ and $\Ab_\YY$. \label{part_fGalois_XcapY}
\end{enumerate}
\end{lemma}
\begin{proof}
%The discussion in \ref{sect_centralalgsuffices}
%showed that $\Ab_\XX$ is a subalgebra of $\Zb(\C[G])$,
%and 
\thmref{lem_subalgyieldsfZs} shows that the map $\XX\mapsto\Ab_\XX$ is invertible.
For part (\ref{part_fGalois_po}),
suppose $\XX\pleq\YY$ and let $Y$ be a part of $\YY$.
Then $Y$ is the union of some parts of $\XX$,
so $f_Y$ is the sum of the corresponding idempotents $f_X$.
Then $f_Y\in\Ab_\XX$; hence $\Ab_\YY\sseq \Ab_\XX$.
Conversely, suppose $\Ab_\YY\sseq \Ab_\XX$.
Let $Y$ be a part of $\YY$.
Then $f_Y$ is an idempotent in $\Ab_\XX$,
so it is a sum of some of the spanning idempotents $\{f_X: X\in\XX\}$ of $\Ab_\XX$.
It follows that $Y$ must be a union of parts of $\XX$, as desired.
Thus the map $\XX\mapsto \Ab_\XX$ is an order-reversing bijection
between $\Part(\Irr(G))$ and the partially ordered set of
the subalgebras of $\Zb(\C[G])$ under inclusion.

Then for part (\ref{part_fGalois_XcupY}),
because $\XX\join\YY$ is the least upper bound for $\XX$ and $\YY$ in $\Part(\Irr(G))$,
the subalgebra $\Ab_{\XX\join\YY}$ must be the largest subalgebra
contained in both $\Ab_\XX$ and $\Ab_\YY$, namely $\Ab_\XX \cap \Ab_\YY.$
Similarly for part (\ref{part_fGalois_XcapY}),
we observe that $\XX\meet\YY$ is the greatest lower bound for $\XX$ and $\YY$,
so $\Ab_{\XX\meet\YY}$ is the smallest subalgebra 
containing both $\Ab_\XX$ and $\Ab_\YY$,
which by definition is $\gen{\Ab_\XX, \Ab_\YY}$.\qed
\end{proof}

The map from the partitions $\KK$ of $G$ 
to the subspaces $\textspanof{\Khat: K\in\KK}$
forms no such bijection,
not even when we require $\KK$ to be coarser than the conjugacy classes
and those subspaces actually to be subalgebras.
The proof of \thmref{lem_subalgyieldsfZs}, 
which established the bijection used to prove \thmref{lem_fGalois},
relied heavily on the orthogonality of the $e_\chi$'s, 
but the $\Khat$'s exhibit no such orthogonality.
An analogue of \thmref{part_fGalois_XcupY} does hold, however.

\begin{lemma}\label{lem_KhalfGalois}
Let $\KK,\LL\in\Part(G)$,
and write $\MM=\KK\join\LL$.
Then
\begin{equation}\label{eqn_intersectKspans}
\spanof{\Mhat: M\in\MM} = \spanof{\Khat: K\in\KK} \cap \spanof{\Lhat: L\in\LL}.
\end{equation}
\end{lemma}
\begin{proof}
%$(\sseq)$
Since each part $M\in\MM=\KK\join\LL$ 
is a union of some parts of $\KK$, %\thmref{sorites2} is used here implicitly
the sum $\Mhat$ lies in $\textspanof{\Khat: K\in\KK}$.
Likewise $\Mhat\in\textspanof{\Lhat: L\in\LL}$,
so the left side of (\ref{eqn_intersectKspans}) is contained in the right hand side.

%$(\supseteq)$
On the other hand, each element $d$ on the right side of (\ref{eqn_intersectKspans})
may be written as 
$d = \sum_{K\in\KK} a_K \Khat = \sum_{L\in\LL} b_L \Lhat$ for some coefficients $a_K,b_L\in\C$.
Recall that each element $g\in G$ occurs in exactly one $K$ and in exactly one $L$,
and that $G$ is a basis for $\C[G]$.
Now the function mapping $g$ to the coefficient of $g$ in $d$ is constant on
each $K\in\KK$, and also constant on each $L\in\LL$.
Hence it is constant on each member of $\KK\join\LL=\MM$ by \thmref{fconstonxcupy},
and so $d$ lies in the span of $\{\Mhat: M\in\MM\}$. \qed
\end{proof}

%\vspace{.2in}
These two lemmas allow us to define a binary ``join'' operation on
supercharacter theories of a group $G$:

\begin{proposition}\label{prop_joinisthy}
  Let $G$ be a group.
  Let %, and let 
  $\XX$ and $\YY$ be partitions of $\Irr(G)$
  and $\KK$ and $\LL$ partitions of $G$
  such that $(\XX, \KK)$ and $(\YY,\LL)$ are each supercharacter theories of $G$.
  Then $(\XX\join \YY\, ,\, \KK\join \LL)$ is also a supercharacter theory of $G$,
  which is denoted $(\XX,\KK)\join(\YY,\LL)$.
\end{proposition}
\begin{proof}
Let $\ZZ = \XX\join\YY$ and $\MM=\KK\join\LL$.
%Note first that since $\{1\}\in\KK$ and $\{1\}\in\LL$, it follows that $\{1\}\in\MM$.
%
To show that the functions $\{\sigma_Z: Z\in\ZZ\}$
are constant on the sets $M\in\MM$,
let $Z\in\ZZ$, let $M\in\MM$, and let $g,h\in M$.
Now $Z={\bigcup}_{X\in \II} X$ for some subset $\II\sseq\XX$,
so $\s_Z = \sum_{X\in \II} \s_X$ must be constant on each $K\in\KK$
because $(\XX,\KK)$ is a supercharacter theory.
On the other hand, by symmetry
$\s_Z$ is also constant on each $L\in\LL$.
So by \thmref{fconstonxcupy}, 
it follows that $\s_Z$ is constant on each $M\in\MM$.
It only remains to show that $|\ZZ|=|\MM|$.

Recall that Diaconis and Isaacs showed \cite[Theorem 2.2(b)]{diaconis_isaacs}
that $\left\{f_X: X \in \XX\right\}$ and $\{\Khat: K\in\KK\}$ 
are two different bases for the algebra $\BB((\XX,\KK))$.
Hence
$$\makebox[0pt][c]{$
  \begin{array}{r@{}c@{}l@{}l}
  \spanof{f_Z: Z\in\ZZ} &{}={}& \spanof{f_X: X\in\XX} \cap \spanof{f_Y: Y\in\YY}  
                                  & \mbox{~(by \thmref{part_fGalois_XcupY})} \\
                      &=& \spanof{\Khat: K\in\KK} \cap \spanof{\Lhat: L\in\LL}
                                  & \mbox{~(by \cite[Theorem 2.2(b)]{diaconis_isaacs})} \\
                                  %& \mbox{~(by \cite[Theorem 2.2(b)]{diaconis_isaacs})} \\
                      &=& \spanof{\Mhat: M\in\MM}.
                                  & \mbox{~(by \thmref{lem_KhalfGalois})} \\
  \end{array}
  $}$$
But since both $\{f_Z: Z\in\ZZ\}$ and
$\{\Mhat: M\in\MM\}$ are linearly independent sets over $\C$, 
both $|\ZZ|$ and $|\MM|$ must equal the dimension of the algebra in question,
so $|\ZZ|=|\MM|$ as desired.
We conclude that $(\XX\join\YY\widecomma \KK\join\LL)$ is a supercharacter theory of $G$.\qed
\end{proof}

If both $(\XX,\KK)$ and $(\YY,\LL)$ are supercharacter theories of $G$,
it is in general \emph{not} true that 
$(\XX\meet\YY\, ,\, \KK\meet\LL)$ is again a supercharacter theory.
For example, let $G$ be the symmetric group on six letters,
whose eleven conjugacy classes
may be denoted by their cycle structures.
%as $1, 2^1, 2^2, 2^3, 3^1, 3^1 2^1, 3^2, 4^1, 4^1 2^1, 5^1, 6^1$.
%
Then the two 
partitions 
$\KK=\left\{ 1,~ 2^1\cup 2^3\cup 4^1,~ 6^1\cup 3^1 2^1,~ 2^2,~ 3^1\cup 3^2,~ 4^1 2^1,~5^1\right\}$
and
$\LL=\left\{ 1,~ 2^1\cup 2^3\cup 4^1\cup 6^1\cup 3^1 2^1,~ 2^2,~ 3^1,~ 3^2,~4^1 2^1,~5^1\right\}$
both correspond to supercharacter theories, but their meet
%$$\KK\meet\LL = \left\{ 1,~ 2^1\cup 2^3\cup 4^1,~ 6^1\cup 3^1 2^1,~ 2^2,~ 3^1,~ 3^2,~4^1 2^1,~5^1\right\}$$
does not.

%\vspace{.2in}
We conclude this discussion of joins 
by making $\Sup(G)$ into a partially ordered set.
Because every supercharacter theory in $\Sup(G)$ 
corresponds to a partition $\XX\in\Part(\Irr(G))$,
one possibility would be to declare that
$(\XX,\KK)\pleq (\YY,\LL)$ if $\XX\pleq\YY$.
On the other hand,
we could just as well partially order $\Sup(G)$ with respect 
to the partition of $G$ given by each theory's superclasses;
according to this method, we would 
write $(\XX,\KK)\pleq(\YY,\LL)$ if $\KK\pleq\LL$.
Thanks to \thmref{prop_joinisthy},
we can show that each of these two alternatives 
produces the same partial ordering of $\Sup(G)$.

\begin{corollary}\label{prop_XlatticeisKlattice}
Let $(\XX,\KK)$ and $(\YY,\LL)$ be supercharacter theories of $G$.
Then $\XX\pleq \YY$ if and only if $\KK\pleq\LL$.
\end{corollary}
\begin{proof}
$(\Leftarrow)$
Suppose $\KK\pleq\LL$.
Since $(\XX,\KK)$ and $(\YY,\LL)$ are supercharacter
theories for $G$, so is 
$(\XX\join\YY \widecomma \KK\join\LL)$, which is equal to $(\XX\join\YY \widecomma \LL)$.
But because superclasses and supercharacters determine one another, %by \thmref{fact_KdeterminesX} 
there is only one supercharacter theory with superclasses $\LL$,
namely $(\YY,\LL)$.  So
$\XX\join\YY = \YY$ and $\XX\pleq\YY$ as desired.

$(\Rightarrow)$:
Suppose $\XX\pleq\YY$.
Then $(\XX\join\YY \widecomma \KK\join\LL) = (\YY\widecomma \KK\join\LL)$
is a supercharacter theory for $G$.
Because $(\YY,\LL)$ is the unique supercharacter theory with supercharacters from $\YY$,
%By \thmref{fact_XdeterminesK},
we must have $\KK\join\LL=\LL$
and $\KK\pleq\LL$.
\qed
\end{proof}

We are therefore not breaking symmetry between superclasses and supercharacters
when we define a partial ordering of $\Sup(G)$ as follows. 

\begin{definition}
  Let $(\XX,\KK)$ and $(\YY,\LL)$ be supercharacter theories of a group $G$.
  Then we write $(\XX,\KK)\pleq(\YY,\LL)$ if $\XX\pleq\YY$.
\end{definition}

\ifdatestamp
  \pagestyle{fancy}
\fi

%\chapter{3}{The $\mathbf{*}$-product}\label{ch_starproduct}

\section{The $\mathbf{*}$-product}\label{sect_starproduct}

In this section, we show that if $N$ is a  normal subgroup of $G$,
then some supercharacter theories of $N$ 
can be combined with
supercharacter theories of $G/N$
to form supercharacter theories of the full group $G$.
We can thus construct supercharacter theories of large groups 
by combining those of smaller groups.

Let a group $G$ act on another group $H$;
then $G/\Cb_G(H)$ embeds naturally in $\Aut(H)$,
so there exists a supercharacter theory
$\Conj{G}{H}=(\XX,\KK)\in\Sup(H)$ 
such that $\XX$ is the partition of $\Irr(H)$ into $G$-orbits
and $\KK$ is the finest partition of $H$ into unions of conjugacy classes
such that each part is $G$-invariant.
Then for another supercharacter theory $(\YY,\LL)\in\Sup(H)$,
every part $L\in\LL$ is $G$-invariant if and only if $\KK\pleq\LL$,
which is true if and only if $\XX\pleq\YY$ (by \thmref{prop_XlatticeisKlattice}), 
which is true
if and only if every part $Y\in\YY$ is $G$-invariant.
We may thus unambiguously speak of a $G$-invariant supercharacter theory.

\begin{definition}
  Let $G$ and $H$ be groups and let $G$ act on $H$ by automorphisms.
  We say that $(\XX,\KK)\in\Sup(H)$ is
  \hbox{\defnstyle{$\mathbf{G}$-invariant}}
  if the action of $G$ fixes each part $K\in\KK$ setwise.
  We denote by $\Supinv{G}(H)$ 
  the set of all $G$-invariant supercharacter theories of $H$.
\end{definition}

For example, if $N\norm G$,
then $\CCC\in\Sup(N)$ is $G$-invariant
if and only if its superclasses are unions of $G$-conjugacy classes.
Also, every supercharacter theory of $G$ is $G$-invariant, 
so $\Supinv{G}(G)=\Sup(G)$.
Note that if $M,N\norm G$ with $N<M$,
then a supercharacter theory of $M/N$ is $G/N$-invariant if and only if it is $G$-invariant.
As our notation indicates,
the supercharacter theory $\Conj{G}{H}$ 
is the minimal $G$-invariant member of $\Sup(H)$.

We would like to define a product
$$*: \Supinv{G}(N)\by\Sup(G/N)\longrightarrow \Sup(G).$$
So suppose $\CCC=(\XX,\KK)\in\Supinv{G}(N)$ and $\DDD=(\YY,\LL)\in\Sup(G/N)$.
Let us first consider the superclasses: 
$\KK$ is a $G$-invariant partition of $N$
and $\LL$ a $G$-invariant partition of $G/N$, 
one part of which is the coset $N$.
We can ``inflate'' $\LL$ to be a partition of $G$,
one part of which will be the set $N$,
which we can then replace with the partition $\KK$ of $N$.
%The resulting partition of $G$ incorporates the information from both $\KK$ and $\LL$,
%and each part is a union of conjugacy classes 
%because both $\KK$ and $\LL$ are $G$-invariant.
To express this formally, 
%we define the following notation:
%
for each subset $L\sseq G/N$ let $\tw{L}$ denote the inflation ${\bigcup}_{Ng\in L}Ng$.
Extend this notation to a set $\LL$ of subsets of $G/N$
by $\tw{\LL}=\{\tw{L}: L\in\LL\}$,
and let $\LL^\circ$ denote $\LL\setminus\{\{N\}\}$.
Then for $\KK\in\Part(N)$ and $\LL\in\Part(G/N)$, we have a partition
\begin{equation}\label{eqn_classpart}
  \KK\cup\tw{\LL^\circ}\in\Part(G).
\end{equation}

%For the supercharacters,
%Continuing to suppose $(\XX,\KK)\in\Supinv{G}(N)$
%and $(\YY,\LL)\in\Sup(G/N)$,
%we want to construct a partition of $\Irr(G)$
%out of $\XX$ and $\YY$ somehow.
%We can view each irreducible character of $G/N$
%as a character of $G$ by inflation, 
%it is natural to let $\YY$ be a subset of our partition.
%It remains to use each set $X\in\XX$ of irreducible characters of $N$
%to construct a set of irreducible characters of $G$.

For the supercharacters, 
recall that if $N$ is a normal subgroup of $G$ and $\psi\in\Irr(N)$,
then $\Irr(G|\psi)$ denotes the set of irreducible characters
$\chi$ of $G$ such that $[\chi_N,\psi]>0$.
If $Z\sseq\Irr(N)$ is a union of $G$-orbits,
then define the subset $Z^G$ of $\Irr(G)$
to be
${\bigcup}_{\psi\in Z} \Irr(G|\psi)$.
Extend this notation to a set $\ZZ$ of subsets of $\Irr(N)$
by letting $\ZZ^G=\{Z^G: Z\in\ZZ\}$,
and let $\ZZ^\circ=\ZZ\setminus\{\{1_N\}\}$.

Now consider $(\XX,\KK)\in\Sup_G(N)$ and $(\YY,\LL)\in\Sup(G/N)$ as before.
Since $\XX$ is a partition of $\Irr(N)$ into unions of $G$-orbits,
it follows that $\XX^G$ is a partition of $\Irr(G)$.
Since $\{1_N\}\in\XX$, one part of $\XX^G$ is 
$\{1_N\}^G=\{\chi\in\Irr(G) : N\sseq\ker\chi\}$, 
which we identify with $\Irr(G/N)$ in the usual way.
Thus we can replace that part of $\XX^G$
with the partition $\YY$ of $\Irr(G/N)$,
obtaining a partition of $G$
\begin{equation}\label{eqn_charpart}
  \YY\cup(\XX^\circ)^G\in\Part(\Irr(G))
\end{equation}
incorporating information from both $\XX$ and $\YY$.

%Note that if $X$ and $Y$ are disjoint,
%then so are $X^G$ and $Y^G$, 
%since we require $X$ and $Y$ to be unions of $G$-orbits.

We shall show that the partitions of (\ref{eqn_classpart}) and (\ref{eqn_charpart})
do form a supercharacter theory of $G$,
by way of a brief \thmtype{lem_calculatesigmaXG} 
demonstrating the suitability of the notation ``$X^G$.''

\begin{lemma}\label{lem_calculatesigmaXG}
Let $N\norm G$ and let $X\sseq \Irr(N)$
be a union of $G$-orbits.
Then 
$$\sigma_{(X^G)} = \left(\sigma_X\right)^G.$$
\end{lemma}
\begin{proof}
  It suffices to prove the statement when $X$ is a single $G$-orbit.
  Let $\ZZ$ be the partition of $\Irr(N)$ into $G$-orbits,
  so that $\ZZ^G$ is a partition of $\Irr(G)$.
  Then the regular character 
  %$\rho_N=\sum_{\mathclap{Z\in\ZZ}} \s_Z$ 
  $\rho_N=\sum_{{Z\in\ZZ}} \s_Z$ 
  and so
  $$\rho_G = (\rho_N)^G=\sum_{Z\in\ZZ} (\s_Z)^G.$$
  Now the characters $(\s_Z)^G$ have no irreducible constituents in common
  with one another, 
  so $(\s_X)^G=\s_Y$ where $Y$ is the set of irreducible constituents of $(\s_X)^G$.
  But $Y=X^G$, so we conclude that
  $\s_{(X^G)} = (\s_X)^G$.\qed
\end{proof}

%Using this notation,
%if $\XX$ is a $G$-invariant partition of $\Irr(N)$ and
%$\YY\in\Part(\Irr(G/N))$, then
%\begin{equation}
%  \YY\cup \left(\XX^\circ\right)^G\in\Part(\Irr(G)).
%\end{equation}

\begin{theorem}\label{thm_starproduct}
  Let $G$ be a group and let $N\norm G$.
  Suppose $(\XX,\KK)\in\Supinv{G}(N)$ and suppose $(\YY,\LL)\in\Sup(G/N)$.
  Then
  $$\left(
       \YY\cup \left(\XX^\circ\right)^G,~
       \KK\cup\tw{\LL^\circ}
    \right)$$
  is a supercharacter theory of $G$.
\end{theorem}
\begin{proof}
Let $\ZZ=\YY\cup \left(\XX^\circ\right)^G\in\Part(\Irr(G))$
and let $\MM=\KK\cup\tw{\LL^\circ}\in\Part(G)$.
Now
$|\ZZ|=|\YY|+|\XX|-1 = |\LL|+|\KK|-1 = |\MM|$,
so it remains to show that for each part $Z\in\ZZ$,
the character $\sigma_Z$ is constant on each subset $M\in\MM$.

One possibility is that $Z$ lies in $\YY$.
In this case, because $Z\sseq\Irr(G/N)$,
the character $\sigma_Z$ has $N$ in its kernel,
so it is certainly constant on every part $K\in\KK$.
Moreover, because $(\YY,\LL)$ is a supercharacter theory of $G/N$,
we know that $\sigma_Z$ (viewed as a character of $G/N$)
is constant on each superclass $L\in\LL$;
viewed as a character of $G$, 
it is therefore constant on each $\tw{L}\in\tw{\LL^\circ}$.
So $\sigma_Z$ is constant on each set $M\in\MM$ in the case that $Z\in\YY$.

The other possibility is that $Z=X^G$ for some part $X\in\XX$.
Now because $(\XX,\KK)$ is $G$-invariant,
we know that $X$ is a union
of $G$-orbits of $\Irr(N)$,
so we can calculate by \thmref{lem_calculatesigmaXG}
that $\sigma_Z = \sigma_{(X^G)}=(\sigma_X)^G$.
Because $N\norm G$,
we know that $(\s_X)^G$ vanishes outside of $N$;
therefore $\s_Z$ is constant on every set $\tw{L}\in\tw{\LL^\circ}$.
Moreover,
because the set $X$ is $G$-invariant, 
the character $\s_X$ of $N$ is invariant under the action of $G$,
so the restriction $(\s_Z)_N=((\s_X)^G)_N=\ix{G}{N} \s_X$.
Because $\s_X$ is constant on every part $K\in\KK$,
so too is $\s_Z$.

Thus $\s_Z$ is constant on $M$ for every $Z\in\ZZ$ and every $M\in\MM$,
so we conclude that $(\ZZ,\MM)$ is a supercharacter theory of $G$.\qed
\end{proof}

We may therefore define a product of supercharacter theories as follows.

\begin{definition}\label{defn_starproduct}
  Let $G$ be a group and let $N\norm G$.
  Let $\CCC=(\XX,\KK)\in\Supinv{G}(N)$ and $\DDD=(\YY,\LL)\in\Sup(G/N)$.
  Then their \defnstyle{$\mathbf{*}$-product}, written $\CCC\star[N]\DDD$,
  is the supercharacter theory of $G$
  $$\left( \YY\cup\left(\XX^\circ\right)^G \, , \, \KK\cup\tw{\LL^\circ} \right).$$
  We say that $\CCC\star[N]\DDD$ is a $*$-product \defnstyle{over} $N$
  or that $\CCC\star[N]\DDD$ \defnstyle{factors over} $N$.
  When $N$ is clear from context, we omit the subscript and write simply $\CCC*\DDD$.
\end{definition}

\section{Factoring}

Let us investigate some properties of the $*$-product.
First, we would like to recognize which supercharacter
theories of a group $G$ arise as $*$-products.
One characteristic of a $*$-product over $N$
is that $N$ is a union of some of the superclasses.
%Let us give this condition a name.

\begin{definition}\label{defn_Cnormal}
  Let $G$ be a group and let $\CCC\in\Sup(G)$.
  Then a subgroup $N$ of $G$ which is a union of some superclasses of $\CCC$
  is called \defnstyle{$\CCC$-normal}.
\end{definition}

Recall that $\m(H)$ denotes the minimal supercharacter theory of a group $H$,
namely its ordinary character theory; 
likewise $\M(H)$ denotes the supercharacter theory with
exactly two superclasses, $1$ and $H-1$.
If $N\norm G$,
then $\M(N)*\M(G/N)$ is a supercharacter theory of $G$
with the three superclasses $1$, $N-1$, and $G-N$;
the corresponding partition of $\Irr(G)$
is $\{\{1_G\},\Irr(G/N)-\{1_G\},\Irr(G)\setminus\Irr(G/N)\}$.
We shall denote this supercharacter theory by $\MMt{N}{G}$.
Let $\CCC=(\XX,\KK)\in\Sup(G)$;
then $N$ is $\CCC$-normal 
if and only if $\CCC\pleq\MMt{N}{G}$,
which is true if and only 
$\Irr(G/N)$ is the union of some members of $\XX$.\label{disc_altCnormalX}
%the partition of $\Irr(G)$ associated with $\CCC$.
%if no supercharacter $\s_X$ associated with $\CCC$ has
%irreducible constituents both inside and outside $\Irr(G/N)$.

A ``minimal'' counterpart to $\MMt{N}{G}$ %$\M(N)*\M(G/N)$
is $\Conj{G}{N}*\m(G/N)$,
whose superclasses are those conjugacy classes of $G$ which lie in $N$,
together with the nontrivial conjugacy classes of $G/N$ pulled back to $G$.
Note that within the partition lattice of $G$,
this is the finest partition of $G$ into unions of conjugacy classes
such that each nontrivial $N$-coset
lies entirely within some superclass.
We shall denote this supercharacter theory by the symbol $\mm{N}{G}$.
In its partition of $\Irr(G)$, 
the characters in $\Irr(G/N)$ are in singleton parts,
while every part outside $\Irr(G/N)$
is of the form $\Irr(G|\psi)$ for some $\psi\in\Irr(N)$.

We noted above that 
if $\CCC$ is a $*$-product over a normal subgroup $N$ of $G$,
then $N$ is \mbox{$\CCC$-normal}.
The converse is not true:
a supercharacter theory $\CCC$
may not factor over $N$
even if $N$ is $\CCC$-normal.
We nevertheless can construct
related supercharacter theories of $N$ and of $G/N$
from such a supercharacter theory.
Since $N$ is $\CCC$-normal,
those superclasses of $\CCC$ which lie in $N$ partition $N$,
and we shall show that this partition belongs to a supercharacter theory $\CCC_N$ of $N$.
Likewise we shall prove that 
those supercharacters of $\CCC$ which have $N$ in their kernels
are the supercharacters of a supercharacter theory $\CCC^{G/N}$ of $G/N$.

To prove these statements,
we need a little notation.
If $Z\sseq\Irr(G)$ is a union of sets of the form $\Irr(G|\psi)$
for various $\psi\in\Irr(N)$,
let $f(Z)$ denote the set of all 
irreducible constituents of $\sum_{\chi\in Z} \chi_N$,
so that $(f(Z))^G=Z$.
Moreover, let $\phi:G\to G/N$ be the canonical homomorphism.

\begin{definition}\label{defn_restrictdownandup}
  Let $G$ be a group, let $\CCC\in\Sup(G)$,
  and let $N$ be a $\CCC$-normal subgroup of $G$.
  Writing $\CCC\join\mm{N}{G}=(\ZZ,\MM)$
  and defining $\phi$ and $f$ as above, 
  let
  $$\CCC_N = 
    \big( \{f(Z): Z\in\ZZ,~Z\not\sseq\Irr(G/N)\}\cup\{\{1_N\}\}
           \widecomma
           \{M\in\MM: M\sseq N\} \big)$$
  and
  $$\CCC^{G/N} = 
    \big( \{Z\in\ZZ: Z\sseq\Irr(G/N)\}
           \widecomma
           \{\phi(M): M\in\MM,~M\not\sseq N\}\cup\{\{N\}\} \big).$$
\end{definition}

Note that replacing $\CCC$ with $\CCC\join\mm{N}{G}$ does not change 
  the superclasses which lie within $N$,
  nor that portion of the partition of $\Irr(G)$ which lies within $\Irr(G/N)$.
Therefore the supercharacters of $\CCC^{G/N}$ are indeed those supercharacters of $\CCC$
with $N$ in their kernels,
and the superclasses of $\CCC_N$ are those superclasses of $\CCC$ which lie in $N$.
In the proof of the following lemma, let us say ``$\XX$ is constant on $\KK$'' to mean
that for each $X\in\XX$, the character $\s_X$ is constant on every part $K\in\KK$.

\begin{lemma}\label{lem_restrictworks}
  Let $N$ be a subgroup of a group $G$, 
  let $\CCC\in\Sup(G)$, and suppose $N$ is $\CCC$-normal.
  Then $\CCC_N$ is a $G$-invariant supercharacter theory of $N$
  and $\CCC^{G/N}$ is a supercharacter theory of $G/N$.
  Moreover, $$\CCC\join\mm{N}{G}=\CCC_N*\CCC^{G/N}.$$
\end{lemma}
\begin{proof}
  Let $\BBB=\CCC\join\mm{N}{G}$ 
  and write $\BBB=(\ZZ,\MM)$.
  Now since both $\CCC\pleq\MMt{N}{G}$ and $\mm{N}{G}\pleq\MMt{N}{G}$,
  it follows that $\BBB=\CCC\join\mm{N}{G}\pleq\MMt{N}{G}$,
  so $N$ is $\BBB$-normal.
  Let
  \begin{eqnarray*}
    \KK &=& \{M\in\MM: M\sseq N\},                                 \\
    \LL &=& \{\phi(M): M\in\MM,~M\not\sseq N\}\cup\{\{N\}\}        \\ 
    \XX &=& \{f(Z): Z\in\ZZ,~Z\not\sseq\Irr(G/N)\}\cup\{\{1_N\}\}, \mbox{~and}\\
    \YY &=& \{Z\in\ZZ: Z\sseq\Irr(G/N)\},                          
  \end{eqnarray*}
  so that by definition $\CCC_N = (\XX,\KK)$ and $\CCC^{G/N}=(\YY,\LL)$.
  Let us now verify that the sets $\XX$, $\KK$, $\YY$, and $\LL$
  are partitions of the appropriate sets.
  
  Since $N$ is $\BBB$-normal, it follows that $\KK=\{M\in\MM: M\sseq N\}$
  is a partition of $N$,
  and the set $\{M\in\MM: M\not\sseq N\}$ is a partition of $G-N$.
  Since $\mm{N}{G}\pleq\BBB$, each element of this latter set is a union of $N$-cosets,
  and thus $\LL^\circ = \{\phi(M): M\in\MM,~M\not\sseq N\}$ is a partition of $G/N-\{N\}$,
  so $\LL$ is a partition of $G/N$.
  Because $\tw{\phi(M)}=M$, we may also note that $\MM=\KK\dcup \tw{\LL^\circ}$.
  
  As for the characters,
  $\YY=\{Z\in\ZZ: Z\sseq\Irr(G/N)$
  is a partition of $\Irr(G/N)$
  because $N$ is $\BBB$-normal.
  Likewise the set $\{Z\in\ZZ: Z\not\sseq\Irr(G/N)\}$ is a partition of $\Irr(G)-\Irr(G/N)$,
  and since $\mm{N}{G}\pleq \BBB$, each $Z$ in this set is a union of sets of the form $\Irr(G|\psi)$.
  Thus $\XX^\circ = \{f(Z): Z\in\ZZ,~Z\not\sseq\Irr(G/N)\}$ is a well-defined partition of $\Irr(N)-\{1_N\}$,
  so $\XX$ is a partition of $\Irr(N)$.
  Moreover, since $(f(Z))^G=Z$, we also see that $\ZZ = \YY\dcup (\XX^\circ)^G$.

  \gobble{
                         Let $\BBB=\CCC\join\mm{N}{G}$ 
                         and write $\BBB=(\ZZ,\MM)$ as in \thmref{defn_restrictdownandup}.
                         Now since both $\CCC\pleq\MMt{N}{G}$ and $\mm{N}{G}\pleq\MMt{N}{G}$,
                         it follows that $\BBB=\CCC\join\mm{N}{G}\pleq\MMt{N}{G}$,
                         so $N$ is $\BBB$-normal.
                         
                         Let $\KK=\{M\in\MM: M\sseq N\}$
                         \marginpar{Do I really need these two paragraphs?}
                         and let $\LL_0=\{\phi(M):M\in\MM,~M\not\sseq N\}$.
                         Because $N$ is $\BBB$-normal, 
                         we know that $\KK$ is a partition of $N$ containing $\{1\}$;
                         since $\mm{N}{G}\pleq\BBB$, 
                         we also know that every nontrivial $N$-coset 
                         is contained in some member of $\MM$,
                         so $\LL_0$ is a partition of $G/N\setminus \{N\}$.
                         We may then form a partition of $G/N$
                         by letting $\LL=\LL_0\cup\{\{N\}\}$;
                         note that $\MM=\KK\dcup\tw{\LL^\circ}$.
                         
                         Now let $\YY=\{Z\in\ZZ: Z\sseq\Irr(G/N)\}$
                         and let $\XX_0=\{f(Z): Z\in\ZZ,~Z\not\sseq\Irr(G/N)\}$.
                         Because $N$ is $\BBB$-normal, we know that $\YY\in\Part(G/N)$.
                         Since $\mm{N}{G}\pleq\BBB$,
                         we also know that the set $\Irr(G|\psi)$ 
                         is contained in some member of $\ZZ$
                         for every nonprincipal $\psi\in\Irr(N)$;
                         hence the members of $\XX_0$ are well-defined and pairwise disjoint,
                         and their union is $\Irr(N)-\{1_N\}$.
                         So let $\XX=\XX_0\cup\{\{1_N\}\}$,
                         which is a partition of $\Irr(N)$.
                         Note that $\ZZ=\YY\dcup (\XX^\circ)^G$.
  }
  
  To show that $\CCC_N$ and $\CCC^{G/N}$ are indeed supercharacter theories,
  it remains to show that the purported supercharacters are actually constant on the superclasses,
  that $|\XX|=|\KK|$, and that $|\YY|=|\LL|$.
  
  Consider first $\CCC_N$. %$(\XX,\KK)$.
  Let $X\in\XX$.  
  If $X=\{1_N\}$, then $\s_X=1_N$ is trivially constant on all $K\in\KK$;
  otherwise, $X=f(Z)$ for some $Z\in\ZZ$,
  so $Z=X^G$.
  Then because
  $$(\s_Z)_N = \left(\s_{(X^G)}\right)_N = \left( (\s_X)^G\right)_N = \ix{G}{N} \s_X$$
  %$$\ix{G}{N} \s_X = \left( (\s_X)^G\right)_N=\left(\s_{X^G}\right)_N = (\s_Z)_N$$
  is constant on each $K\in\KK$, so too is $\s_X$.
  We conclude that $\XX$ is constant on $\KK$.
  
  Now consider $(\YY,\LL)$,
  and let $Y\in\YY\sseq\ZZ$.  
  Then $\s_Y$ is constant 
  on every superclass $M\in\MM$,
  and in particular, on those superclasses outside $N$.
  But since $\s_Y$ has $N$ in its kernel,
  when viewed as a character of $G/N$
  it is constant on the images of those superclasses,
  namely the members of $\LL^\circ$.
  Moreover, $\s_Y$ is certainly constant on the singleton set $\{N\}\sseq G/N$.
  We conclude that $\YY$ is constant on $\LL$.
  
  Now because $\XX$ is constant on $\KK$
  and $\YY$ is constant on $\LL$,
  by \cite[Theorem 2.2]{diaconis_isaacs}
  we know that $|\XX|\leq |\KK|$ and $|\YY|\leq |\LL|$.
  But then
  $$|\KK|+|\YY|-1\leq |\KK|+|\LL|-1=|\MM|=|\ZZ|=|\YY|+|\XX|-1\leq |\YY|+|\KK|-1,$$
  so equality must hold throughout; hence $|\XX|=|\KK|$ and $|\YY|=|\LL|$.
  %But we also have that
  %$$|\YY|+|\XX|-1=|\ZZ|=|\MM|=|\KK|+|\LL|-1,$$
  %so $|\XX|+|\YY| = |\KK|+|\LL|$.
  %Hence equality must hold in both cases,
  %and $|\XX|=|\KK|$ and $|\YY|=|\LL|$.
  We conclude that $\CCC_N=(\XX,\KK)$ is a supercharacter theory of $N$
  and $\CCC^{G/N}=(\YY,\LL)$ a supercharacter theory of $G/N$;
  the former is $G$-invariant 
  because its superclasses are also superclasses of $\CCC$.
  Finally, by definition
  $$\CCC_N*\CCC^{G/N}
        = (\XX,\KK)*(\YY,\LL)
        = \left(\YY\cup (\XX^\circ)^G \widecomma \KK \cup \tw{\LL^\circ} \right)
        = (\ZZ, \MM)
        = \CCC\join\mm{N}{G}$$
  as desired.\qed
\end{proof}

With the help of the preceding \thmtype{lem_restrictworks},
we can determine whether a supercharacter theory $\EEE$ of $G$
factors over a normal subgroup $N$.

\begin{corollary}\label{cor_factoringcriteria}
  Let $G$ be a group, let $N\norm G$, and
  let $\EEE$ be a supercharacter theory of $G$.
  Then $\EEE$ factors over $N$ if and only if
    $N$ is $\EEE$-normal and
    every superclass outside $N$ is a union of $N$-cosets.
  Moreover, if $\EEE=\CCC \star[N] \DDD$,
  then $\CCC=\EEE_N$ and $\DDD=\EEE^{G/N}$.
\end{corollary}
\begin{proof}
  Suppose $\EEE=\CCC\star[N]\DDD$;
  write $\CCC=(\XX,\KK)$ and $\DDD=(\YY,\LL)$,
  so that by definition
  $\EEE=\left( \YY\cup (\XX^\circ)^G \widecomma \KK\cup\tw{\LL^{\circ}}\right)$.
  Then $N={\bigcup}_{{K\in\KK}} K$ is $\EEE$-normal,
  and each superclass of $\EEE$ outside $N$ 
  is the preimage of some $L\in\LL$,
  and hence a union of $N$-cosets.
  Moreover, the superclasses of $\EEE_N$
  are those superclasses of $\EEE$ that lie in $N$, namely the superclasses of $\CCC$;
  thus $\EEE_N=\CCC$.
  Likewise the supercharacters of $\EEE^{G/N}$ 
  are the supercharacters of $\EEE$ with $N$ in their kernels, 
  namely the supercharacters of $\DDD$,
  so $\EEE^{G/N}=\DDD$.
  
  Now suppose for the converse that 
  $N$ is $\EEE$-normal and 
  that every superclass of $\EEE$ outside $N$ is a union of $N$-cosets.
  Then $\mm{N}{G}\pleq\EEE$, so
  $$\EEE = \EEE\join\mm{N}{G}=\EEE_N * \EEE^{G/N}$$
  by \thmref{lem_restrictworks}.\qed
\end{proof}

Since $\EEE\pleq\MMt{N}{G}$ if and only if $N$ is $\EEE$-normal
and $\mm{N}{G}\pleq\EEE$ if and only if every superclass outside $N$
is a union of $N$-cosets,
the preceding \thmtype{cor_factoringcriteria}
can be rephrased in terms of the partial order on $\Sup(G)$:
$\EEE$ factors over $N$
if and only if 
$\mm{N}{G}\pleq \EEE \pleq \MMt{N}{G}$.

\section{Associativity}

Up to this point, 
we have been working with a fixed normal subgroup $N$ of $G$.
Suppose now that we have two normal subgroups $N$ and $M$,
with $1\leq N\leq M\leq G$,
and three supercharacter theories $\CCC\in\Sup(N)$, $\DDD\in\Sup(M/N)$, and $\EEE\in\Sup(G/M)$.
We may summarize this situation in a diagram:
%$$ 1 \underbrace{\leq}_{\CCC} N \underbrace{\leq}_{\DDD} M \underbrace{\leq}_{\EEE} G.$$
\newlength{\halfN}\newlength{\halfM}
\newlength{\myleni}\newlength{\mylenii}\newlength{\myleniii}
\settowidth{\halfN}{$N$}\setlength{\halfN}{0.5\halfN}
\settowidth{\halfM}{$M$}\setlength{\halfM}{0.5\halfM}
\settowidth{\myleni}{$1{~}\leq {~}N$} \addtolength{\myleni}{-\halfN} \addtolength{\myleni}{-.3ex}
\settowidth{\mylenii}{$N {~}\leq {~} M$} 
     \addtolength{\mylenii}{-\halfN}
     \addtolength{\mylenii}{-\halfM} 
     \addtolength{\mylenii}{-.3ex}
\settowidth{\myleniii}{$M{~}\leq G$}
     \addtolength{\myleniii}{-\halfM}
     \addtolength{\myleniii}{-.3ex}
$$ \makebox[0pt][l]{$1 {~} \leq {~} N{~} \leq {~} M {~} \leq G.$}\hspace{-1pt}
   \underbrace{\rule{\myleni}{0pt}}_{\CCC}
   \underbrace{\rule{\mylenii}{0pt}}_{\DDD}
   \underbrace{\rule{\myleniii}{0pt}}_{\EEE}
$$
Suppose $\CCC$ and $\DDD$ are $G$-invariant.
Then in particular $\CCC$ is $M$-invariant, 
so we can form the product $\CCC\star[N]\DDD\in\Sup(M)$,
which is also $G$-invariant.
Thus  we can form $(\CCC\star[N]\DDD)\star[M]\EEE\in\Sup(G)$.
On the other hand, 
we can also form $\DDD\star[M/N]\EEE\in\Sup(G/N)$
and then form $\CCC\star[N](\DDD\star[M/N]\EEE)\in\Sup(G)$ as well.
Fortunately, the two products $(\CCC * \DDD)*\EEE$ and $\CCC*(\DDD*\EEE)$ are the same; 
in other words, the $*$-product is associative.

\begin{lemma}\label{lem_starassoc}
Let $G$ be a group and let $N$ and $M$ be normal subgroups with $N\leq M$. %\leq G$.
Suppose $\CCC\in\Supinv{G}(N)$, $\DDD\in\Supinv{G}(M/N)$, and $\EEE\in\Sup(G/M)$.
Then
$$(\CCC \star[N] \DDD) \star[M] \EEE = \CCC \star[N] (\DDD \star[M/N] \EEE).$$
\end{lemma}
\begin{proof}
  Let $\CCC = (\WW,\JJ)$, let $\DDD = (\XX,\KK)$, and let $\EEE = (\YY,\LL).$
  Consider the partitions of characters.
  For $\CCC* \DDD$, the partition of characters is $\XX\cup(\WW^\circ)^M$;
  thus for $(\CCC*\DDD)*\EEE$,
  the partition of characters is
  $$   \YY \cup \left( \left( \XX\cup (\WW^\circ)^M\right)^\circ\right)^G
     = \YY \cup \left( \XX^\circ \cup (\WW^\circ)^M \right)^G
     = \YY \cup (\XX^\circ)^G \cup (\WW^\circ)^G.
     $$
  On the other hand,
  the partition of characters for $\DDD*\EEE$
  is $\YY\cup (\XX^\circ)^G$,
  so the partition of characters for $\CCC*(\DDD*\EEE)$ is
  also $\YY\cup(\XX^\circ)^G \cup (\WW^\circ)^G$.
  Thus $(\CCC \star[N] \DDD) \star[M] \EEE$ and $\CCC \star[N] (\DDD \star[M/N] \EEE)$
  are two supercharacter theories of $G$ with identical supercharacters, and hence are equal.\trueqed
\end{proof}

%
%
%\begin{proof}
%  %Our setup is the following:
%  %$$ 1 \underbrace{\leq}_{\CCC} N \underbrace{\leq}_{\DDD} M \underbrace{\leq}_{\EEE} G.$$
%  Let $\CCC = (\WW,\JJ)$, let $\DDD = (\XX,\KK)$, and let $\EEE = (\YY,\LL).$
%  %
%  For this proof only, we must adjust the notation of \thmref{defn_Lpullback}
%  to differentiate between two operations:
%  if $S$ is a subset of $G/N$ or $G/M$,
%  then the symbol $\tw{S}$ will denote the preimage of $S$ in $G$ as before,
%  but the symbol $\vartw{T}$ will denote the preimage in $G/N$ of a set $T\sseq G/M$.
%  %(Note that $\tw{\vartw{T}}=\tw{T}$.)
%  %
%  Now the superclasses of $\CCC\star[N]\DDD$ are 
%  $\JJ\cup \tw{\KK^\circ}$,
%  so the superclasses of $(\CCC\star[N]\DDD)\star[M]\EEE$ are
%  $$\left( \JJ \cup \tw{\KK^\circ} \right) \cup \tw{\LL^\circ}
%    = \JJ \cup \tw{\KK^\circ} \cup \tw{\LL^\circ}.$$
%
%  On the other hand, the superclasses of $\DDD\star[M/N]\EEE$ are
%  $\KK\cup\vartw{\LL^\circ}$,
%  and then the superclasses of $\CCC\star[N](\DDD\star[M/N]\EEE)$ are
%  $$\JJ\cup \tw{\left(\KK\cup\vartw{\LL^\circ}\right)^\circ} 
%    = \JJ\cup \tw{\left(\KK^\circ \cup\vartw{\LL^\circ}\right)} 
%    = \JJ\cup \left(\tw{\KK^\circ} \cup\tw{\vartw{\LL^\circ}}\right)
%    = \JJ\cup \tw{\KK^\circ} \cup\tw{\LL^\circ}
%  $$
%  Hence $(\CCC\star[N]\DDD)\star[M]\EEE$ 
%  and $\CCC\star[N](\DDD\star[M/N]\EEE)$
%  are supercharacter theories 
%  of $G$ with identical superclasses, so they are equal.\qed
%\end{proof}

\section{Unique factorization}

Although \thmref{cor_factoringcriteria} showed that a supercharacter theory $\EEE$ of $G$
can be written as a $*$-product over $N$ either in a unique way or not at all,
\thmref{lem_starassoc} shows that a supercharacter theory might factor over two
different normal subgroups $N$ and $M$.
However, this can only happen if $N$ contains $M$ or vice versa.
%To show this, we first need a lemma.

\begin{lemma}\label{cor_MNcontainment}
Let $G$ be a group with normal subgroups $N$ and $M$.
Let $\CCC\in\Sup(G)$,
and suppose $\CCC$ is a $*$-product both over $N$ and over $M$.
Then either $N\leq M$ or $M\leq N$.
\end{lemma}
\begin{proof}
Suppose $N\not\leq M$, and choose an element $n\in N\setminus M$.
Then because $\CCC$ factors over $M$ and $n\not\in M$, 
the entire coset $Mn$ must lie in the same superclass of $\CCC$ as $n$.
But then because $N$ is $\CCC$-normal,
we must have $Mn\sseq N$ and thus $M\leq N$.\qed
\end{proof}

Let us investigate further what happens when $\CCC\star[N]\DDD=\EEE\star[M]\FFF$.
Without loss of generality, 
we have normal subgroups $N$ and $M$ of $G$
with $1\leq N\leq M\leq G$,
where $\CCC\in\Supinv{G}(N)$, $\DDD\in\Sup(G/N)$, \mbox{$\EEE\in\Supinv{G}(M)$},
and $\FFF\in\Sup(G/M)$.
We may portray this setup in a diagram:
\begingroup
%\newlength{\halfN}\newlength{\halfM}
\newlength{\mylenC}\newlength{\mylenD}\newlength{\mylenE}\newlength{\mylenF}
\settowidth{\halfN}{$N$}\setlength{\halfN}{0.5\halfN}
\settowidth{\halfM}{$M$}\setlength{\halfM}{0.5\halfM}
\settowidth{\mylenC}{$1{~}\leq {~}N$}
     \addtolength{\mylenC}{-\halfN} 
     \addtolength{\mylenC}{-.0ex}
\settowidth{\mylenD}{$N {~}\leq {~} M{~}\leq G$} 
     \addtolength{\mylenD}{-\halfN}
     \addtolength{\mylenD}{-.0ex}
\settowidth{\mylenE}{$1{~}\leq {~}N{~}\leq {~}M$}
     \addtolength{\mylenE}{-\halfM}
     \addtolength{\mylenE}{-.3ex}
\settowidth{\mylenF}{$M{~}\leq {~}G$}
     \addtolength{\mylenF}{-\halfM}
     \addtolength{\mylenF}{-.3ex}
$$
\super{%
       \overbrace{\rule{\mylenC}{0pt}}^{\CCC}%
       \overbrace{\rule{\mylenD}{0pt}}^{\DDD}%
}%
\raisebox{-0.5ex}{\makebox[0pt][l]{$\underbrace{\rule{\mylenE}{0pt}}_{\EEE}
                  \underbrace{\rule{\mylenF}{0pt}}_{\FFF}
                  $}}
1{~}\leq {~}N{~} \leq {~}M{~} \leq {~}G
$$
\endgroup
In this situation, we shall show that 
there exists a $G$-invariant supercharacter theory
\mbox{$\GGG\in\Supinv{G}(M/N)$} 
such that our supercharacter theory of $G$ is simply $\CCC\star[N]\GGG\star[M]\FFF$,
with $\EEE=\CCC\star[N]\GGG$
and $\DDD=\GGG\star[M/N]\FFF$,
as in this diagram:
\begingroup
\newlength{\mylenG}
\settowidth{\halfN}{$N$}\setlength{\halfN}{0.5\halfN}
\settowidth{\halfM}{$M$}\setlength{\halfM}{0.5\halfM}
\settowidth{\mylenC}{$1{~}\leq {~}N$}
     \addtolength{\mylenC}{-\halfN} 
     \addtolength{\mylenC}{-.3ex}
\settowidth{\mylenE}{$1{~}\leq {~}N{~}\leq {~}M$}
     \addtolength{\mylenE}{-\halfM}
     \addtolength{\mylenE}{-.3ex}
\settowidth{\mylenD}{$N {~}\leq {~} M{~}\leq G$} 
     \addtolength{\mylenD}{-\halfN}
     \addtolength{\mylenD}{-.0ex}
\settowidth{\mylenF}{$M{~}\leq {~}G$}
     \addtolength{\mylenF}{-\halfM}
     \addtolength{\mylenF}{-.3ex}
\settowidth{\mylenG}{$N{~}\leq {~}M$}
     \addtolength{\mylenG}{-\halfN}
     \addtolength{\mylenG}{-\halfM}
     \addtolength{\mylenG}{-.3ex}
$$
\super{\overbrace{\rule{\mylenE}{0pt}}^{\EEE}}%
\super[2.9ex]{\rule{\mylenC}{0pt} \rule{2pt}{0pt}%
       \overbrace{\rule{\mylenD}{0pt}}^{\DDD}}%
\raisebox{-0.5ex}{\makebox[0pt][l]{$%
                  \underbrace{\rule{\mylenC}{0pt}}_{\CCC}
                  \underbrace{\rule{\mylenG}{0pt}}_{\GGG}
                  \underbrace{\rule{\mylenF}{0pt}}_{\FFF}
                  $}}
1 {~}\leq{~}N{~} \leq {~} M {~} \leq {~} G.
$$
\endgroup

\begin{lemma}\label{lem_uniquefactorizationlemma}
Let $G$ be a group, let $M,N\norm G$,
and suppose $\CCC \star[N] \DDD = \EEE \star[M] \FFF\in\Sup(G)$.
Without loss of generality, suppose $N\leq M$.
Then there exists some supercharacter theory $\GGG\in\Supinv{G}(M/N)$ 
such that $\EEE = \CCC \star[N] \GGG$
and $\DDD=\GGG \star[M/N] \FFF$.
\end{lemma}
\begin{proof}
First we shall show that $\EEE$ factors over $N$.
Let $\BBB=\CCC \star[N] \DDD = \EEE \star[M] \FFF$.
Now the superclasses of $\EEE$ are 
precisely those superclasses of $\BBB$
which lie in $M$.
Because $\BBB$ factors over $N$, we know that $N$ is $\BBB$-normal;
hence $N$ is $\EEE$-normal.
Also because $\BBB$ factors over $N$,
every superclass of $\BBB$ outside $N$ is a union of $N$-cosets.
In particular, every superclass of $\EEE$ outside $N$ is a union of $N$-cosets.
Then by \thmref{cor_factoringcriteria},
we know $\EEE$ must be a $*$-product over $N$.
Because the superclasses of $\EEE$ lying in $N$ are exactly
the superclasses of $\BBB$ lying in $N$, namely the superclasses of $\CCC$,
we have $\EEE=\CCC\star[N]\GGG$ for some supercharacter theory $\GGG\in\Sup(M/N)$.
Since $\CCC$ is $G$-invariant, so too is $\GGG$.

Then $\BBB = \EEE\star[M]\FFF = (\CCC \star[N] \GGG) \star[M] \FFF
           = \CCC \star[N] (\GGG \star[M/N] \FFF)$
by \thmref{lem_starassoc};
but we also know $\BBB = \CCC \star[N] \DDD$.  
By \thmref{cor_factoringcriteria}, 
we conclude that $\DDD = \GGG\star[M/N]\FFF$, as desired.\qed
\end{proof}

The important implication of these lemmas is that 
every supercharacter theory of a group $G$ 
can be factored uniquely into
a $*$-product of one or more supercharacter theories 
that cannot themselves be written as $*$-products
in a nontrivial way.

\begin{definition}
Let $G$ be a group and let $\CCC\in\Sup(G)$.
We say $\CCC$ is \defnstyle{decomposable}
if $\CCC$ is a $*$-product over a proper nontrivial normal subgroup of $G$.
We say $\CCC$ is \defnstyle{indecomposable}
if it is not decomposable and $|G|>1$.
\end{definition}

\begin{theorem}\label{thm_uniquefactorization}
Let $G$ be a nontrivial group and let $\CCC\in\Sup(G)$.
Then there exists a unique chain of normal subgroups
$1=N_0<N_1<\cdots<N_r=G$ 
and unique indecomposable supercharacter theories
$\DDD_i\in\Supinv{G}(N_i/N_{i-1})$ for $i=1,\ldots,r$
such that
$$\CCC=\DDD_1 *\DDD_2*\cdots*\DDD_r.$$
\end{theorem}
\begin{proof}
Among all chains of normal subgroups $1=N_0<N_1<\cdots<N_r=G$
for which there exist (possibly decomposable) supercharacter theories 
$\DDD_i\in\Supinv{G}(N_i/N_{i-1})$
such that $\CCC=\DDD_1*\cdots*\DDD_r$,
choose one of maximal length $r$.  
This can be done because $G$ is finite.

Now if $\DDD_j$ were decomposable for some $j\in\{1,\ldots,r\}$,
then there would exist
a normal subgroup $M$ of $N_j$ such that $N_{j-1}<M<N_j$,
with the property that 
$\DDD_j=\EEE\star[M/N_{j-1}]\FFF$ 
for some $\EEE\in\Sup_{N_j}(M/N_{j-1})$ and $\FFF\in\Sup(N_j/M)$.
Since $\DDD_j$ is $G$-invariant, 
in fact $M$ would be normal in $G$ and both $\EEE$ and $\FFF$ would be $G$-invariant.
Then 
$$1=N_0<\cdots<N_{j-1}<M<N_j<\cdots<N_r=G$$
would be a chain of length $r+1$ with
$\CCC=\DDD_1*\cdots*\DDD_{j-1}*\EEE*\FFF*\DDD_{j+1}*\cdots*\DDD_r$, 
contradicting the maximality of our choice.
Hence each $\DDD_i$ is indecomposable,
and we have proven the existence of a factorization
into indecomposable supercharacter theories.

Now to show uniqueness, induct on $|G|$.
Suppose there are two chains of normal subgroups 
$1=N_0<\cdots<N_r=G$ and $1=M_0<\cdots<M_s=G$
and two sets of indecomposable supercharacter theories
$\DDD_i\in\Sup(N_i/N_{i-1})$ for $i=1,\ldots,r$ 
and $\EEE_i\in\Sup(M_i/M_{i-1})$ for $i=1,\ldots,s$
such that 
$$\CCC=\DDD_1*\cdots*\DDD_r=\EEE_1*\cdots*\EEE_s.$$

Suppose $N_1\neq M_1$; then 
without loss of generality, by \thmref{cor_MNcontainment}
we may assume that $N_1<M_1$.
Then we have
\begingroup
\newlength{\halfNone}
\settowidth{\halfNone}{$N_1$}\setlength{\halfNone}{0.5\halfNone}
\newlength{\halfMone}
\settowidth{\halfMone}{$M_1$}\setlength{\halfMone}{0.5\halfMone}
\newlength{\mylenDone}
\settowidth{\mylenDone}{$1 {~}< {~} N_1$} 
     \addtolength{\mylenDone}{-\halfNone}
     \addtolength{\mylenDone}{-.2ex}
\newlength{\mylenDrest}
\settowidth{\mylenDrest}{$N_1{~}<{~}M_1{~}\leq{~}G$}
     \addtolength{\mylenDrest}{-\halfNone}
     \addtolength{\mylenDrest}{-.0ex}
\newlength{\mylenEone}
\settowidth{\mylenEone}{$1 {~}< {~} N_1 {~}<{~}M_1$} 
     \addtolength{\mylenEone}{-\halfMone}
     \addtolength{\mylenEone}{-.2ex}
\newlength{\mylenErest}
\settowidth{\mylenErest}{$M_1{~}\leq{~}G$}
     \addtolength{\mylenErest}{-\halfMone}
     \addtolength{\mylenErest}{-.0ex}
$$
\super{\overbrace{\rule{\mylenDone}{0pt}}^{\DDD_1}
       \overbrace{\rule{\mylenDrest}{0pt}}^{\DDD_2*\cdots*\DDD_r}
}%
\raisebox{-0.5ex}{\makebox[0pt][l]{$%
                  \underbrace{\rule{\mylenEone}{0pt}}_{\EEE_1}
                  \underbrace{\rule{\mylenErest}{0pt}}_{\EEE_2*\cdots*\EEE_s}
                  $}}
1{~}<{~}N_1{~}<{~}M_1{~}\leq{~}G
$$
\endgroup
and by \thmref{lem_uniquefactorizationlemma} we know $\EEE_1$ must factor over $N_1$,
contradicting the indecomposability of $\EEE_1$.
Hence $N_1=M_1$, so
$$\CCC=\DDD_1 \star[N_1] \left(\DDD_2*\cdots*\DDD_r\right)
  \mbox{~and~}
  \CCC=\EEE_1 \star[N_1] \left(\EEE_2*\cdots*\EEE_s\right)$$
are identical products over $N_1$.
Then \thmref{cor_factoringcriteria} implies 
that $\DDD_1=\EEE_1$ and that $\DDD_2*\cdots*\DDD_r=\EEE_2*\cdots*\EEE_s$.
By applying the inductive hypothesis to $G/N$,
we see both that $r-1=s-1$ 
and that $N_i/N_1=M_i/N_1$ and $\DDD_i=\EEE_i$ for all $i\in\{2,\ldots,r\}$.
Therefore $r=s$ and $N_i=M_i$ and $\DDD_i=\EEE_i$ for all $i\in\{1,\ldots,r\}$,
proving the uniqueness of the factorization.\qed
\end{proof}
 
\ifdatestamp
  \pagestyle{fancy}
\fi

%\chapter{4}{The $\bigwtp$-product}\label{ch_wtp}

\section{The $\mathbf{\bigwtp}$-product}\label{sect_wtp}

We have investigated the $*$-product
constructed from a supercharacter theory $\CCC$ 
of a normal subgroup $M\norm G$
and a supercharacter theory $\DDD$ of the quotient group $G/M$.
A similar construction can still be done in the more general situation
when $\DDD$ is a supercharacter theory of a quotient of $G$ by a \emph{smaller} normal subgroup $N$,
provided that $\CCC$ and $\DDD$ satisfy certain conditions.
We are considering the situation of the following diagram:
$$\makebox[0pt][l]{$\underbrace{\phantom{1\leq N\leq M}}_{\CCC}$}1\leq \overbrace{N\leq M\leq G}^{\DDD}$$
In order to put $\CCC$ and $\DDD$ together to form a supercharacter theory for $G$,
we will of course want $N$ to be $\CCC$-normal and $M/N$ to be $\DDD$-normal,
but we will also want the ``overlap'' of the two theories on $M/N$ to be the same;
more explicitly, we will require $\CCC^{M/N}=\DDD_{M/N}$.

\begin{theorem}\label{thm_wtpproduct}
  Let $G$ be a group with subgroups $N\leq M\leq G$.
  Suppose $\CCC\in\Sup_G(M)$ and $\DDD\in\Sup(G/N)$ such that
  \begin{enumerate}
    \item $N$ is $\CCC$-normal,
    \item $M/N$ is $\DDD$-normal, and
    \item $\CCC^{M/N}=\DDD_{M/N}$.
  \end{enumerate}
  Then there exists a unique supercharacter theory $\EEE\in\Sup(G)$
  such that $\EEE_M=\CCC$ and $\EEE^{G/N}=\DDD$
  and every superclass outside $M$ is a union of $N$-cosets.

  Using our earlier notation,
  if $\CCC=(\XX,\KK)$ and $\DDD=(\YY,\LL)$, 
  then
  \begin{equation*}%\label{eqn_wtpdefn}
    \EEE=\left( \YY\cup\{X^G: X\in\XX,~X\not\sseq\Irr(M/N)\}
                \widecomma 
                \KK\cup\{\tw{L}: L\in\LL,~L\not\sseq M/N\}\right).
  \end{equation*}
\end{theorem}
\begin{proof}
  For every superclass $L$ of $\DDD$ lying outside $M/N$, 
  take its preimage $\tw{L}$ in $G$;
  because $M/N$ is $\DDD$-normal,
  this gives a partition of $G\setminus M$.
  To this set add all the superclasses of $\CCC$;
  since these partition $M$,
  the resulting set $\KK\cup \{\tw{L}: L\in\LL,~L\not\sseq M/N\}$
  is a partition of $G$
  %with $\{1\}$ as one of the parts,
  which we shall call $\JJ$.
  Recalling that $|\CCC|$ denotes the number of superclasses of $\CCC$,
  note that $|\JJ|=|\CCC|+\left(|\DDD|-|\DDD_{M/N}|\right)$.
  
  Now because $N$ is $\CCC$-normal,
  the subset $\Irr(M/N)$ is a union of parts of $\XX$,
  as discussed in Section \ref{disc_altCnormalX}.
  Hence $\{X\in\XX: X\not\sseq\Irr(M/N)\}$ partitions $\Irr(M)-\Irr(M/N)$,
  so the set $\{X^G: X\in\XX,~X\not\sseq\Irr(M/N)\}$ partitions $\Irr(G)-\Irr(G/N)$
  since $\CCC$ is $G$-invariant.
  Since $\YY$ is a partition of $\Irr(G/N)$,
  the union $\YY\cup\{X^G: X\in\XX,~X\not\sseq\Irr(M/N)\}$
  is a partition of $\Irr(G)$; call it $\WW$.
  %(Symbolically, $\WW=\YY\cup\{X^G: X\in\XX,~X\not\sseq\Irr(M/N)\}$.)
  Note that 
  $$|\WW|=|\DDD|+(|\CCC|-|\CCC^{M/N}|)=|\CCC|+\left(|\DDD|-|\DDD_{M/N}|\right)=|\JJ|.$$
  Then to prove that $(\WW,\JJ)$ is a supercharacter theory of $G$,
  it remains only to show that $\s_W$ is constant on $J$ 
  for each $W\in\WW$ and each $J\in\JJ$.
  
  Let $W\in\WW$.  It may be that $W\in\YY$, 
  so that $\s_W$ is a supercharacter of $\DDD$.
  In this case, there are three sorts of sets $J\in\JJ$ to consider:
  those that lie within $G\setminus M$, those within $M\setminus N$, and those within $N$.
  First, the supercharacter $\s_W$ of $\DDD$
  is constant on each superclass $L$ of $\DDD$ lying outside $M/N$,
  so $\s_W$ (viewed as a character of $G$) 
  is constant on each preimage $\tw{L}$ in $G\setminus M$.
  Thus $\s_W$ is constant on each set $J\in\JJ$ that lies in $G\setminus M$.
  Next note that
  $\s_W$ is constant on the nontrivial superclasses of $\DDD_{M/N}=\CCC^{M/N}$.
  Therefore $\s_W$ is constant on the preimages of these superclasses,
  which are exactly the superclasses of $\CCC_N*\CCC^{M/N}$ which lie outside $N$.
  Now $\CCC\pleq \CCC\join\mm{N}{G} = \CCC_N*\CCC^{M/N}$,
  %but since $\CCC_N*\CCC^{M/N}=\CCC\join\mm{N}{G}$,
  %those preimages are exactly the superclasses of $\CCC\join\mm{N}{G}$ outside $N$.
  so every superclass $K$ of $\CCC$ lying outside of $N$
  is contained within a superclass of $\CCC_N*\CCC^{M/N}$ outside $N$;
  thus $\s_W$ is constant on that superclass $K$.
  Hence $\s_W$ is constant on every set $J\in\JJ$ that lies in $M\setminus N$.
  Finally, $\s_W$ has $N$ in its kernel, 
  so it is constant on those sets $J\in\JJ$ that lie within $N$.
  Therefore $\s_W$ is constant on every member of $\JJ$,
  under the supposition that $W\in\YY$.
  
  The other possibility is that $W=X^G$ 
  for some part $X\in\XX$ not lying in $\Irr(M/N)$.
  Since $\CCC$ is $G$-invariant, the set $X$ must be a union of $G$-orbits,
  so $\s_W=\s_{(X^G)}=(\s_X)^G$ by \thmref{lem_calculatesigmaXG}.
  Then because $M\norm G$, we know that $\s_W$ vanishes outside $M$,
  and hence is constant on all parts $J\in\JJ$ that lie outside $M$.
  On the other hand, each part $J\in\JJ$ that lies in $M$ 
  is a superclass of $\CCC$,
  and when $\s_W$ is restricted to $M$, 
  the character $(\s_W)_M=\left((\s_X)^G\right)_M=\ix{G}{M}\s_X$ is constant
  on $J$ because $\s_X$ is.
  
  Hence $\s_W$ is constant on each part $J\in\JJ$ for all parts $W\in\WW$, 
  and we conclude that $(\WW,\JJ)$ is a supercharacter theory of $G$.  
  Let $\EEE=(\WW,\JJ)$; we need to show 
  that $\EEE$ satisfies the conclusions of the \thmtype{thm_wtpproduct}.
  By construction, the superclasses of $\EEE$ that lie in $M$
  are the superclasses of $\CCC$, so $\CCC=\EEE_M$.
  Likewise the supercharacters of $\EEE^{G/N}$
  are those supercharacters of $\EEE$ that have $N$ in their kernels,
  namely the supercharacters of $\DDD$;
  hence $\DDD=\EEE^{G/N}$.
  Third, by construction
  the superclasses of $\EEE$ outside $M$ are preimages
  of certain superclasses of $\DDD$, so they are unions of $N$-cosets.
  
  Finally, to show uniqueness, suppose $\FFF\in\Sup(G)$
  %Finally, to show uniqueness, suppose a supercharacter theory $\FFF\in\Sup(G)$
  satisfies the conditions that \mbox{$\FFF_M=\CCC$},
  that $\FFF^{G/N}=\DDD$,
  and that every superclass of $\FFF$ outside $M$
  is a union of $N$-cosets.
  %We shall show that $\FFF=\EEE$ by comparing their superclasses.
  Then $\FFF_M=\CCC=\EEE_M$, 
  so $\EEE$ has the same superclasses within $M$ as does $\FFF$.
  Moreover, because the superclasses of $\FFF$ outside of $M$ are unions
  of $N$-cosets, the set
  $$\begin{array}{r@{}c@{}l}
    \{\mbox{superclasses of $\FFF$ outside $M$}\}
      &{}={}& \{\mbox{superclasses of $\FFF\join\mm{N}{G}$ outside $M$}\} \\
      &=& \{\mbox{superclasses of $\FFF_N*\FFF^{G/N}$ outside $M$}\} \\
      &=& \{\mbox{superclasses of $\FFF_N*\DDD$ outside $M$}\} \\
      &=& \{\mbox{preimages of the superclasses of $\DDD$ outside $M/N$}\} \\
      &=& \{\mbox{superclasses of $\EEE$ outside $M$}\}.
  \end{array}$$
  Therefore $\FFF$ has the same superclasses as $\EEE$,
  so $\FFF=\EEE$ as desired.\qed
\end{proof}

We may therefore define a $\bigwtp$-product as follows.

\begin{definition}\label{defn_triangleproduct}
  Let $G$ be a group with subgroups $N$ and $M$
  such that $N\leq M$,
  and suppose $\CCC=(\XX,\KK)\in\Sup_G(M)$ and $\DDD=(\YY,\LL)\in\Sup(G/N)$.
  If $N$ is $\CCC$-normal, if $M/N$ is $\DDD$-normal,
  and if $\CCC^{M/N}=\DDD_{M/N}$,
  then we define their \defnstyle{$\bigwtp$-product},
  written $\CCC\wtp\DDD$, to be the supercharacter theory of $G$
  $$\left( \YY\cup\{X^G: X\in\XX,~X\not\sseq\Irr(M/N)\}
                \widecomma 
                \KK\cup\{\tw{L}: L\in\LL,~L\not\sseq M/N\}\right).$$
  We say that $\CCC\wtp\DDD$ is a \defnstyle{$\bigwtp$-product over $N$ and $M$}.
  If $1<N$ and $M<G$, then we say that $\CCC\wtp\DDD$
  is a \defnstyle{nontrivial $\bigwtp$-product}.
\end{definition}

%\section{Recognition}

If $N=M$, then the condition that $\CCC^{M/N}=\DDD_{M/N}$ 
holds trivially
and $\CCC\wtp\DDD=\CCC*\DDD$.  
Thus the $\bigwtp$-product is a generalization of the $*$-product.
An analogue of \thmref{cor_factoringcriteria} holds for $\bigwtp$-products,
giving necessary and sufficient conditions for a supercharacter theory
to be a $\bigwtp$-product.

\begin{proposition}\label{prop_wtprecognition}
  Let $G$ be a group, let $\EEE\in\Sup(G)$,
  and let $N$ and $M$ be $\EEE$-normal subgroups of $G$
  with $N\leq M$.
  Then $\EEE$ is a $\bigwtp$-product over $N$ and $M$
  if and only if every superclass outside $M$ is a union of $N$-cosets.
  In this case, $\EEE=\EEE_M\wtp\EEE^{G/N}$.
\end{proposition}
\begin{proof}
  Suppose $\EEE=\CCC\wtp\DDD$ over $N$ and $M$.
  Then by definition
  the superclasses of $\EEE$ outside $M$ are unions of $N$-cosets.
  
  So suppose for the converse 
  that every superclass of $\EEE$ outside $M$ 
  is a union of \mbox{$N$-cosets}.
  Let $\CCC=\EEE_M$ and let $\DDD=\EEE^{G/N}$;
  we want to show that $\EEE=\CCC\wtp\DDD$,
  but first we need to show that this $\bigwtp$-product is defined.
  Because $N$ and $M$ are $\EEE$-normal,
  it follows that $N$ is $\CCC$-normal and $M/N$ is $\DDD$-normal.
  %Recall that $\EEE\join \mm{N}{G} = \EEE_N*\EEE^{G/N}=\EEE_N*\DDD$.
  %Therefore the nontrivial superclasses of $\DDD$ lying within $M/N$
  %are the images modulo $N$ of the superclasses of $\EEE\join\mm{N}{G}$ in $M\setminus N$.
  %But these are the same as the images of the superclasses 
  %of $\CCC\join\mm{N}{G}$ in $M\setminus N$,
  %which are the superclasses of $\CCC^{M/N}$
  %since $\CCC\join\mm{N}{M}=\CCC_N*\CCC^{M/N}$.
  %Hence $\DDD_{M/N}=\CCC^{M/N}$.
  %
  Let $\phi: G\to G/N$ be the quotient homomorphism.  Then the set
  $$\begin{array}{r@{}c@{}l}
    \{\mbox{nontrivial superclasses of $\DDD_{M/N}$}\}
      &{}={}& \{\mbox{nontrivial superclasses of $\DDD$ in $M/N$}\} \\
      &=& \phi\left( \{\mbox{superclasses of $\EEE_N*\DDD$ in $M\setminus N$}\}\right) \\
      &=& \phi\left( \{\mbox{superclasses of $\EEE_N*\EEE^{G/N}$ in $M\setminus N$}\}\right) \\
      &=& \phi\left( \{\mbox{superclasses of $\EEE\join\mm{N}{G}$ in $M\setminus N$}\}\right) \\
      &=& \phi\left( \{\mbox{superclasses of $\CCC\join\mm{N}{M}$ in $M\setminus N$}\}\right) \\
      &=& \phi\left( \{\mbox{superclasses of $\CCC_N*\CCC^{M/N}$ in $M\setminus N$}\}\right) \\
      &=& \{\mbox{nontrivial superclasses of $\CCC^{M/N}$}\}. \\
  \end{array}$$
  Therefore $\CCC^{M/N}=\DDD_{M/N}$, so we can form the product $\CCC\wtp\DDD$.

  But then $\EEE$ is a supercharacter theory of $G$
  with $\EEE_M=\CCC$ and $\EEE^{G/N}=\DDD$ 
  and every superclass outside $M$ a union of $N$-cosets,
  and \thmref{thm_wtpproduct} guarantees that there is only one such supercharacter theory,
  namely $\CCC\wtp\DDD$.
  So $\EEE=\CCC\wtp\DDD$ as desired.\qed
\end{proof}

\ifdatestamp
  \pagestyle{fancy}
\fi

\section{Dual supercharacter theories}\label{sect_duality}

We conclude this article by 
restricting our attention to abelian groups
and investigating a bijection
between $\Sup(G)$ and $\Sup(\Irr(G))$,
constructed using the
natural isomorphism
$\sim: G\to \Irr(\Irr(G))$.
Recall that $\tw{g}(\chi)$ is defined to be $\chi(g)$ 
for all $\chi\in\Irr(G)$
and all $g\in G$.
If $K$ is a subset of $G$, then we define $\tw{K}$ to be $\{\tw{g}: g\in K\}$;
likewise if $\KK$ is a partition of $G$,
we define $\tw{\KK}$ to be $\{\tw{K}: K\in\KK\}$,
which is a partition of $\Irr(\Irr(G))$.
We shall show that if $(\XX,\KK)\in\Sup(G)$,
then $(\tw{\KK},\XX)\in\Sup(\Irr(G))$.
The proof requires 
us to define a matrix corresponding to a partition of a set.
\begin{definition}
  Let $S$ be a set of size $n$
  and let $\RR$ be a partition of $S$ into $k$ parts.
  Fix an ordering $S=\{s_1,\ldots,s_n\}$
  and an ordering $\RR=\{R_1,\ldots,R_k\}$.
  Then the \defnstyle{partition matrix}
  of $\RR$ is the $k\by n$ matrix $\Rm$ given by 
  $$\Rm_{ij}=\left\{\begin{array}{ll}
                  1, & \mbox{~if~} s_j\in R_i \\
                  0, & \mbox{~if~} s_j\not\in R_i.
             \end{array}\right. $$
\end{definition}

%Note that although the matrix $\Rm$ depends on the orderings chosen,
%it is determined by $\RR$ up to reordering of its rows and columns.
Let $G$ be an abelian group of order $n$
and fix orderings $G=\{g_1,\ldots,g_n\}$ and $\Irr(G)=\{\chi_1,\ldots,\chi_n\}$.
Now $\C[G]$ has two different bases,
an element basis $\{g_1,\ldots,g_n\}$
and an idempotent basis $\{e_{\chi_1},\ldots,e_{\chi_n}\}$.
Therefore there exists a nonsingular $n\by n$ change-of-basis matrix $\Tm$ 
such that if $\xv$ is a row vector 
giving the idempotent coordinates of some element $x$ of $\C[G]$,
then the row vector $\xv\Tm$ gives the element coordinates for $x$.

\begin{lemma}\label{lem_matrixcondition}
  Let $G$ be an abelian group, 
  let $\KK\in\Part(G)$, and let $\XX\in\Part(\Irr(G))$.
  Fix orderings of $G$, $\Irr(G)$, $\KK$, and $\XX$.
  Let $\Km$ and $\Xm$ be the partition matrices of $\KK$ and $\XX$, respectively.
  Let $\Tm$ be the change of basis matrix 
  from the idempotent coordinates to the element coordinates.
  Then $(\XX,\KK)$ is a supercharacter theory of $G$
  if and only if
  $$\rs(\Km)=\rs(\Xm\Tm).$$
\end{lemma}
\begin{proof}
  If $(\XX,\KK)\in\Sup(G)$, 
  then $\textspanof{\Khat: K\in\KK}=\textspanof{f_X: X\in\XX}$.
  On the other hand, if $\textspanof{\Khat: K\in\KK}=\textspanof{f_X: X\in\XX}$,
  then in particular the subspace $\textspanof{\Khat: K\in\KK}$ 
  is a subalgebra of $\Zb(\C[G])$,
  so by \thmref{prop_centralalgsuffices} there is some partition $\YY$ of $\Irr(G)$
  such that $(\YY,\KK)\in\Sup(G)$.  But then
  $$\spanof{f_Y: Y\in\YY} = \spanof{\Khat:K\in\KK} = \spanof{f_X: X\in\XX},$$
  so $\XX=\YY$ by \thmref{lem_subalgyieldsfZs}.
  Thus $(\XX,\KK)$ is a supercharacter theory of $G$ if and only if
  $\textspanof{\Khat: K\in\KK}=\textspanof{f_X: X\in\XX}$.
  
  Now the rows of $\Xm$ are the idempotent coordinates
  of the members of $\{f_X: X\in\XX\}$,
  so the rows of $\Xm\Tm$ are the element coordinates of those same idempotent sums.
  Likewise the rows of $\Km$ are the element coordinates
  of the members of $\{\Khat: K\in\KK\}$.
  So $\{f_X: X\in\XX\}$ and $\{\Khat: K\in\KK\}$ have the same linear span
  if and only if the matrices $\Xm\Tm$ and $\Km$ have the same rowspace.\qed
\end{proof}

With the aid of this \thmtype{lem_matrixcondition},
we now can prove that each supercharacter theory of an abelian group $G$
corresponds to a supercharacter theory of $\Irr(G)$.

\begin{theorem}\label{thm_dualexists}
  Let $G$ be an abelian group,
  and let $(\XX,\KK)$ be a supercharacter theory of $G$.
  Then $(\tw{\KK},\XX)$ is a supercharacter theory of $\Irr(G)$.
\end{theorem}
\begin{proof}
  Let $n=|G|$ and $k=|\KK|$.
  Fix orderings of $G$, $\Irr(G)$, $\XX$, and $\KK$,
  and let $\Xm$ and $\Km$ be the partition matrices
  corresponding to $\XX$ and $\KK$, respectively.
  Order $\tw{G}=\Irr(\Irr(G))$ and $\tw{\KK}$
  in the natural way by letting $\tw{g}_i = \tw{(g_i)}$
  and $\tw{K}_i=\tw{(K_i)}$;
  then the partition matrix corresponding to $\tw{\KK}$
  is also $\Km$.
  
  Now for the algebra $\C[G]$,
  the matrix $\Tm$ which changes the basis
  from idempotent coordinates to element coordinates
  is the $n\by n$ matrix whose $i$th row consists of the element coordinates of 
  $e_{\chi_i}=\rec{|G|}\sum_{g\in G}\bar{\chi_i(g)}\,g$; 
  %$e_{\chi_i}$; 
  hence
  $$\Tm=\left(\rec{|G|} \bar{\chi_i(g_j)}\right)_{ij}
       = \rec{|G|}\bar{ \left(\chi_i(g_j)\right)_{ij}}
       = \rec{|G|}\,\bar{\Cm},$$
  where $\Cm$ is the character table of $G$ viewed as a matrix.
  
  \def\Sm{\mathbf{S}}
  
  On the other hand, the character table of the group $\Irr(G)$
  viewed as a matrix is just $\Cm\transpose$.
  Hence the argument above, when applied to the group algebra $\C[\Irr(G)]$, 
  shows that the matrix $\Sm$ which changes coordinates from the 
  idempotent basis $\{e_{\tw{g}_1},\ldots,e_{\tw{g}_n}\}$
  to the element basis $\{\chi_1,\ldots,\chi_n\}$ is
  $\rec{|G|} \bar{\Cm\transpose}$.
  Since $\bar{\Cm\transpose}=\Cm\inv$
  by the orthogonality of the irreducible characters,
  we have $\Sm= \rec{|G|} \Cm\inv$.

  %Now it happens that
  %\begin{eqnarray*}
  %\rec{|G|} \Cm\Cm\transpose
  %      = \rec{|G|} \big(\chi_i(g_k)\big)_{ik} \big(\chi_j(g_k)\big)_{kj}
  %      &=& \rec{|G|} \left( \sum_{k=1}^n \chi_i(g_k) \chi_j(g_k) \right)_{ij} \\
  %      &=& \left( \rec{|G|} \sum_{g\in G} \chi_i (g) \chi_j(g) \right)_{ij} \\
  %      &=& \big( [\chi_i, \bar{\chi_j}] \big)_{ij} \\
  %      &=& \left( \delta_{\chi_i,\bar{\chi_j}} \right)_{ij},
  %\end{eqnarray*}
  %so multiplication by $\drec{|G|}\Cm\Cm\transpose$
  %permutes the columns of a matrix 
  %by swapping the column corresponding to each irreducible character $\chi$
  %with the column corresponding to $\bar{\chi}$.
  %
  %But since complex conjugation is a field automorphism,
  %interchanging the irreducible characters of $G$ 
  %with their conjugates permutes the parts of $\XX$
  %by \thmref{lem_fieldautspermuteX},
  %so multiplying $\Xm$ by $\drec{|G|}\Cm\Cm\transpose$
  %merely permutes the rows of $\Xm$.

  Now because $(\XX,\KK)\in\Sup(G)$,
  we have $\rs(\Km)=\rs(\Xm\Tm)$
  by \thmref{lem_matrixcondition},
  so $\rs(\Km)=\rs(\Xm\bar{\Cm})$.
  Now both $\Km$ and $\Xm$ are real matrices, 
  so taking complex conjugates of both sides gives
  $\rs(\Km)=\rs(\Xm\Cm)$.
  If two matrices have the same row space,
  then right multiplication by any matrix
  yields two matrices which again have the same row space.
  Therefore 
  $$\rs(\Km\Cm\inv)=\rs(\Xm\Cm\Cm\inv)=\rs(\Xm).$$
  Finally, because $\Sm$ is a scalar multiple of $\Cm\inv$, we have
  $\rs(\Km\Sm)=\rs(\Xm)$.
  %\begin{eqnarray*}
  %  \rs(\Km\Sm)
  %     = \rs(\Km\Cm\inv)
  %     = \rs(\Xm).
  %\end{eqnarray*}
  Then because $\Km$ is the partition matrix of $\tw{\KK}$,
  we conclude by \thmref{lem_matrixcondition}
  that $(\tw{\KK},\XX)$ is a supercharacter theory of $\Irr(G)$.\qed
\end{proof}

\begin{definition}
  Let $G$ be an abelian group
  and let $\CCC=(\XX,\KK)\in\Sup(G)$.
  Then the \defnstyle{dual supercharacter theory}
  $\stw{\CCC}$ is defined to be $(\tw{\KK},\XX)\in\Sup(\Irr(G))$.
\end{definition}

The word ``dual'' is appropriate because
$\stw{\stw{\CCC}}=(\tw{\XX},\tw{\KK})$,
where $\tw{\XX}$ denotes the image of $\XX$
under the natural isomorphism from $\Irr(G)$ to $\Irr(\Irr(\Irr(G)))$.
Thus $\stw{\stw{\CCC}}$ is exactly
the image of $\CCC$ under the natural isomorphism.

\begin{corollary}\label{cor_dualisbij}
  Let $G$ be an abelian group.
  Then the map $\CCC\mapsto\stw{\CCC}$
  defines a bijection from $\Sup(G)$ to $\Sup(\Irr(G))$.
\end{corollary}
\begin{proof}
  Suppose $\CCC,\DDD\in\Sup(G)$ such that $\stw{\CCC}=\stw{\DDD}$;
  then the superclasses of $\stw{\CCC}$ are the same as those of $\stw{\DDD}$.
  Then $\CCC$ and $\DDD$ correspond to the same partition of $\Irr(G)$,
  so $\CCC=\DDD$.  Therefore the map $\CCC\mapsto\stw{\CCC}$ is injective.
  But $|\Sup(G)|=|\Sup(\Irr(G))|$ 
  because $G$ and $\Irr(G)$ are isomorphic,
  so the map must be a bijection.\qed
\end{proof}

Let $(\XX,\KK)$ be a supercharacter theory of an abelian group $G$.
If $X\in\XX$ and $K\in\KK$, then
we know that $\s_X(g)=\sum_{\chi\in X}\chi(g)$
is constant for all $g\in K$.
Thus in the submatrix of the character table 
whose rows lie in $X$ and whose columns lie in $K$,
the column sums are identical.
\thmref{thm_dualexists} shows that 
$\s_{\tw{K}}(\chi)=\sum_{g\in K} \tw{g}(\chi) = \sum_{g\in K}\chi(g)$
is constant for all $\chi\in X$,
so the row sums of that submatrix are also all identical.
Letting $z$ be the sum of all entries in the submatrix,
however, we see that each row sum is ${z}/{|X|}$
and each column sum is ${z}/{|K|}$;
thus the row sums will not equal the column sums
unless $|K|=|X|$.

\section{Duality of $\CCC$-normal subgroups and $\bigwtp$-products}\label{sect_dualprops}

We next investigate the behavior of $\CCC$-normal subgroups
and $\bigwtp$-products under the dual map.
Let $G$ be an abelian group,
and recall that every subgroup of $\Irr(G)$
is of the form $\Irr(G/N)$ for some subgroup $N\leq G$.
If $\CCC\in\Sup(G)$,
a similar connection holds for $\CCC$-normal subgroups.

\begin{lemma}\label{lem_Cnormaldual}
  Let $G$ be an abelian group, let $N$ be a subgroup of $G$,
  and let $\CCC\in\Sup(G)$.
  Then $N$ is $\CCC$-normal if and only if $\Irr(G/N)$ 
  is a $\stw{\CCC}$-normal subgroup of $\Irr(G)$.
\end{lemma}
\begin{proof}
  Write $\CCC=(\XX,\KK)$.
  Recall from Section \ref{disc_altCnormalX}
  that $N$ is $\CCC$-normal if and only if 
  $\Irr(G/N)$ is a union of members of $\XX$,
  which by definition means that $\Irr(G/N)$ is $\stw{\CCC}$-normal.\qed
  %if and only if $\Irr(G/N)$ is a union of parts of $\XX$.\qed
\end{proof}

Moreover, there is a strong connection between 
the supercharacter theories of $G$ 
that arise as $\bigwtp$-products
and the $\bigwtp$-product supercharacter theories of $\Irr(G)$.
We shall show that $\CCC\in\Sup(G)$ is a $\bigwtp$-product
if and only if $\stw{\CCC}$ is a $\bigwtp$-product.
To compute the factors of $\stw{\CCC}$ explicitly, 
we define a new map.
Let $G$ be abelian, let $M\leq G$, and let $\theta\in\Irr(M)$.
Since $G$ is abelian,
the set $\Irr(G|\theta)$ of irreducible characters of $G$ that lie over $\theta$
consists of all extensions of $\theta$ to $G$.
By Gallagher's Theorem \cite[Corollary 6.17]{isaacs_ct},
this set is a coset of $\Irr(G/M)$ in $\Irr(G)$.
Moreover, the map
\begin{equation}\label{eqn_irriso}
  \begin{array}{rccl}
    \irriso: & \Irr(M) & \longrightarrow & \Irr(G)/\Irr(G/M) \\
             & \theta  & \longmapsto     & \Irr(G|\theta)
  \end{array}%
  %\addsymbol{\irriso}{Standard isomorphism of groups from $\Irr(M)$ to $\Irr(G)/\Irr(G/M)$}
\end{equation}
is an isomorphism of groups.
If $X$ is a subset of $\Irr(M)$,
let $X^\irriso=\{\theta^\irriso: \theta\in X\}$;
note that the set $X^G$
is the full preimage of $X^\irriso$ 
with respect to the canonical map from $\Irr(G)\to \Irr(G)/\Irr(G/M)$.
If $\XX$ is a set of subsets of $\Irr(M)$,
let $\XX^\irriso=\{X^\irriso: X\in\XX\}$.

Now if $\stw{\CCC}$ is a supercharacter theory of $\Irr(M)$, 
then the superclasses of $\stw{\CCC}$ partition $\Irr(M)$,
so their images under $\irriso$ partition $\Irr(G)/\Irr(G/M)$.
Because $\irriso$ is a group isomorphism, these images
are the superclasses of a supercharacter theory
of $\Irr(G)/\Irr(G/M)$.

To find the partition of $\Irr(\Irr(G)/\Irr(G/M))$
corresponding to this supercharacter theory,
note that
$\Irr(\Irr(G)/\Irr(G/M))$ 
is the set of all irreducible characters of $\Irr(G)$
that have $\Irr(G/M)$ in the kernel.
But $\tw{g}(\psi)=\psi(g)$ equals $1$ for all $\psi\in\Irr(G/M)$
if and only if $g\in~\bigcap_{{\psi\in\Irr(G/M)}}~\ker\psi = M$.
Thus $\Irr(\Irr(G)/\Irr(G/M))$ is the image of $M$ under the map
$\sim:G\to\Irr(\Irr(G))$.
So if $\KK$ is a partition of $M$,
it follows that $\tw{\KK}$ is a partition of $\Irr(\Irr(G)/\Irr(G/M))$.
We are now ready for the following \thmtype{lem_crosssuperchars}.

\begin{lemma}\label{lem_crosssuperchars}
  Let $G$ be an abelian group, let $M$ be a subgroup of $G$,
  and let $\irriso$ be the isomorphism from $\Irr(M)$ to $\Irr(G)/\Irr(G/M)$ 
  defined in (\ref{eqn_irriso}).
  Let $\CCC=(\XX,\KK)$ be a supercharacter theory of $M$.
  Then $(\tw{\KK},\XX^\irriso)$ is a supercharacter theory of $\Irr(G)/\Irr(G/M)$.
\end{lemma}
\begin{proof}
  We have seen that $\XX^\irriso$ is a partition of $\Irr(G)/\Irr(G/M)$
  because $\XX$ is a partition of $\Irr(M)$.
  The discussion above has also established that 
  $\tw{\KK}$ is a partition of $\Irr(\Irr(G)/\Irr(G/M))$
  because $\KK\in\Part(M)$.
  Certainly $|\tw{\KK}|=|\KK|=|\XX|=|\XX^\irriso|$,
  so it suffices to show that $\s_{\tw{K}}$ is constant on $X^\irriso$
  for each part $K\in\KK$ and each part $X\in\XX$.
  
  Let $K\in\KK$ and let $X\in\XX$.
  For the duration of this proof, 
  let $\alttw{\phantom{C}}$ denote the natural isomorphism from $M$ to $\Irr(\Irr(M))$;
  by \thmref{thm_dualexists}
  $(\alttw{\KK},\XX)$ is a supercharacter theory of $\Irr(M)$,
  so there exists a complex number $c$ 
  such that $\s_{\alttw{K}}(\theta)=c$ for all characters $\theta\in X$.
  Now every member of $X^\irriso$ is 
  %a coset 
  of the form $\chi \Irr(G/M)$
  where $\chi_M\in X$;
  hence
  \begin{eqnarray*}
  \s_{\tw{K}}(\chi\Irr(G/M))
     = \s_{\tw{K}}(\chi)
     = \sum_{m\in K} \tw{m}(\chi)
     = \sum_{m\in K} \chi(m) 
   &=& \sum_{m\in K} \chi_M(m) \\
   &=& \sum_{m\in K} \alttw{m}(\chi_M) \\
   &=& \s_{\alttw{K}}(\chi_M) \\
   &=& c
  \end{eqnarray*}
  is the same for all members $\chi\Irr(G/M)$ of $X^\irriso$.
  Thus $\tw{\KK}$ is constant on $\XX^\irriso$,
  and we conclude that $(\tw{\KK},\XX^\irriso)$ is a supercharacter theory
  of $\Irr(G)/\Irr(G/M)$.\qed
\end{proof}

Let us give a name to the supercharacter theory of \thmref{lem_crosssuperchars}.

\begin{definition}
  Let $G$ be an abelian group, let $M$ be a subgroup of $G$,
  and let $\CCC=(\XX,\KK)$ be a supercharacter theory of $M$,
  so that $\stw{\CCC}=(\tw{\KK},\XX)\in\Sup(\Irr(M))$.
  Then $\stw{\CCC}^\irriso$ denotes the supercharacter theory
  $(\tw{\KK},\XX^\irriso)$ of $\Irr(G)/\Irr(G/M)$.
\end{definition}

Using this definition,
%With this definition in hand,
if a supercharacter theory $\EEE$ of an abelian group is a $\bigwtp$-product,
at last we can write $\stw{\EEE}$ as a $\bigwtp$-product 
and give the factors explicitly.

\begin{proposition}\label{prop_wtpdual}
  Let $G$ be an abelian group
  with subgroups $N\leq M\leq G$,
  and let $\irriso$
  be the isomorphism from $\Irr(M)\to\Irr(G)/\Irr(G/M)$ defined in (\ref{eqn_irriso}).
  Let $\CCC\wtp\DDD$ be a $\bigwtp$-product over $N$ and $M$.
  Then 
  $$\stw{\CCC\wtp\DDD}=\stw{\DDD}\wtp{\stw{\CCC}}^\irriso,$$
  where the second $\bigwtp$-product is over $\Irr(G/M)$ and $\Irr(G/N)$.
\end{proposition}
\begin{proof}
  Let $\CCC=(\XX,\KK)$ and $\DDD=(\YY,\LL)$.
  The relevant subgroups of $\Irr(G)$ are
    $$1 \leq \Irr(G/M) \leq \Irr(G/N) \leq \Irr(G),$$
  and the superclasses of $\stw{\CCC\wtp\DDD}$ are the members of the set
  \begin{equation}\label{eqn_classesCwtpD}
    \YY\cup\{X^G: X\in\XX,~X\not\sseq\Irr(M/N)\}.
  \end{equation}
  Now every superclass outside $\Irr(G/N)$,
  being of the form $X^G$ for some part $X\in\XX$,
  is a union of sets of the form $\Irr(G|\psi)$ where $\psi\in\Irr(M)$,
  and therefore is a union of $\Irr(G/M)$-cosets.
  Then by \thmref{prop_wtprecognition} we know that
  $\stw{\CCC\wtp\DDD}$ factors over $\Irr(G/M)$ and $\Irr(G/N)$ as
  \begin{equation}\label{eqn_dualwtpfactor0}
  \stw{\CCC\wtp\DDD} = \left(\stw{\CCC\wtp\DDD}\right)_{\Irr(G/N)}
                       \wtp
                       \left(\stw{\CCC\wtp\DDD}\right)^{\Irr(G)/\Irr(G/M)}.
  \end{equation}
  But as we saw in (\ref{eqn_classesCwtpD}),
  the superclasses of $\stw{\CCC\wtp\DDD}$ lying in $\Irr(G/N)$
  are the members of $\YY$, which are the superclasses of $\stw{\DDD}$.
  Hence $\left(\stw{\CCC\wtp\DDD}\right)_{\Irr(G/N)}=\stw{\DDD}$.

  On the other hand,
  let $\ZZ$ be the partition of $\Irr(\Irr(G))$ corresponding to $\stw{\CCC\wtp\DDD}$;
  then the members of $\ZZ$ are the images of the superclasses of $\CCC\wtp\DDD$
  under the natural isomorphism $\sim: G \to \Irr(\Irr(G))$.
  We saw in the discussion before \thmref{lem_crosssuperchars}
  that $\Irr(\Irr(G)/\Irr(G/M))$ is exactly $\tw{M}$,
  so the members of $\ZZ$ lying in $\Irr(\Irr(G)/\Irr(G/M))$ are
  the images of those superclasses of $\CCC\wtp\DDD$ that lie in $M$,
  namely the superclasses of $\CCC$.
  Hence the partition of $\Irr(\Irr(G)/\Irr(G/M))$ corresponding to
  $\left(\stw{\CCC\wtp\DDD}\right)^{\Irr(G)/\Irr(G/M)}$
  is exactly $\tw{\KK}$;
  but 
  by \thmref{lem_crosssuperchars}
  this is also the partition of $\Irr(\Irr(G)/\Irr(G/M))$
  corresponding to $\stw{\CCC}^\irriso$.
  Therefore 
  $\left(\stw{\CCC\wtp\DDD}\right)^{\Irr(G)/\Irr(G/M)}=\stw{\CCC}^\irriso$.
  We conclude from Eq.~(\ref{eqn_dualwtpfactor0}) that
  $\stw{\CCC\wtp\DDD} = \stw{\DDD}
                        \wtp
                        \stw{\CCC}^\irriso$,
  as desired.\qed
\end{proof}

\begin{corollary}\label{cor_dualwtpiff}
  Let $G$ be an abelian group with subgroups $N\leq M\leq G$.
  Let $\EEE$ be a supercharacter theory of $G$.
  Then $\EEE$ is a $\bigwtp$-product over $N$ and $M$
  if and only if $\stw{\EEE}$ is a $\bigwtp$-product over $\Irr(G/M)$
  and $\Irr(G/N)$.
\end{corollary}
\begin{proof}
  If $\EEE$ is a $\bigwtp$-product over $N$ and $M$,
  then \thmref{prop_wtpdual} implies that $\stw{\EEE}$ 
  is a $\bigwtp$-product over $\Irr(G/M)$ and $\Irr(G/N)$.
  So now suppose for the converse that $\stw{\EEE}$ 
  is a \mbox{$\bigwtp$-product} over $\Irr(G/M)$ and $\Irr(G/N)$.
  Write $\EEE=(\XX,\KK)$, so that $\stw{\EEE}=(\tw{\KK},\XX)$.
  Then applying \thmref{prop_wtpdual} 
  with $\stw{\EEE}$ in the place of $\EEE$,
  we may conclude that the supercharacter theory 
  $\stw{\stw{\EEE}}=(\tw{\XX},\tw{\KK})$ of $\Irr(\Irr(G))$
  is a $\wtp$-product over $\Irr(\Irr(G)/\Irr(G/N))=\tw{N}$
  and $\Irr(\Irr(G)/\Irr(G/M))=\tw{M}$.
  %(Here $\tw{\XX}$ means the image of $\XX$ under the
  %natural isomorphism $\sim: \Irr(G)\to\Irr(\Irr(\Irr(G)))$.)

  But then the members of $\tw{\KK}$ are the superclasses
  of a $\bigwtp$-product over $\tw{N}$ and $\tw{M}$.
  Since the map $\sim: G\to\Irr(\Irr(G))$ is a group isomorphism,
  we conclude that $\KK$ is the set of superclasses
  of a $\bigwtp$-product over $N$ and $M$.
  Hence $\EEE$ is a $\bigwtp$-product over $N$ and $M$, 
  completing the proof.\qed
\end{proof}

The most important application of \thmref{prop_wtpdual} and \thmref{cor_dualwtpiff}
occurs when $N=M$, in which case we have the following immediate corollary:
\begin{corollary}\label{cor_starproductdual}
  Let $G$ be an abelian group with a subgroup $M\leq G$,
  and let $\irriso$
  be the isomorphism from $\Irr(M)$ to $\Irr(G)/\Irr(G/M)$ defined in (\ref{eqn_irriso}).
  Let $\EEE\in\Sup(G)$;
  then $\EEE$ is a $*$-product over $M$ 
  if and only if $\stw{\EEE}$ is a $*$-product over $\Irr(G/M)$.
  Moreover,
  if $\CCC\in\Sup(M)$ and $\DDD\in\Sup(G/M)$,
  then 
  $$\stw{\CCC\star[M]\DDD}=\stw{\DDD}\star[\Irr(G/M)]{\stw{\CCC}}^\irriso.$$
\end{corollary} 
%\begin{proof}
%  This is exactly %the content of 
%  \thmref{prop_wtpdual} and \thmref{cor_dualwtpiff}
%  under the additional supposition that $N=M$.\qed
%\end{proof}

\section{Conclusion}

We have presented five operations for constructing new supercharacter
theories out of existing ones:
the direct product ($\by$),
the join operation ($\join$),
the $*$-product and its generalization the $\wtp$-product,
and the dual operation $\stw{\phantom{c}}$.
In forthcoming papers we shall show that these operators,
together with the original supercharacter theory constructions
given by Diaconis and Isaacs in \cite{diaconis_isaacs},
suffice to produce all the supercharacter theories
of certain infinite families of finite groups,
including cyclic $p$-groups of odd order
and cyclic groups of order $pq$ and $pqr$.

\bibliography{construction_bib}

\end{document}